\documentclass[11pt]{article}

\usepackage[utf8]{inputenc}
\usepackage[T1]{fontenc}
\usepackage{amsmath,amssymb,amsthm,mathtools}
\usepackage[margin=2.5cm]{geometry}
\usepackage{graphicx}
\usepackage{placeins}
\usepackage{booktabs,longtable,array}
\usepackage{enumitem}
\usepackage{xcolor}
\usepackage[affil-it]{authblk}
\usepackage[colorlinks=true,linkcolor=blue!50!black,citecolor=green!40!black,urlcolor=blue!50!black]{hyperref}

\newtheorem{theorem}{Theorem}[section]
\newtheorem{lemma}[theorem]{Lemma}
\newtheorem{proposition}[theorem]{Proposition}
\newtheorem{corollary}[theorem]{Corollary}
\newtheorem{mainthm}{Theorem}

\newtheorem{maincor}[mainthm]{Corollary}
\theoremstyle{definition}
\newtheorem{definition}[theorem]{Definition}

\theoremstyle{remark}
\newtheorem{remark}[theorem]{Remark}

\newcommand{\R}{\mathbb{R}}
\newcommand{\Z}{\mathbb{Z}}

\newcommand{\NN}{\mathcal{N}}
\newcommand{\D}{D}
\newcommand{\Nm}{\operatorname{N}}
\newcommand{\abs}[1]{\left|#1\right|}
\newcommand{\norm}[1]{\left\|#1\right\|}
\newcommand{\card}{\#}
\newcommand{\Sq}{\mathcal{S}}
\newcommand{\Hc}{\mathcal{H}}
\newcommand{\X}{\mathcal{X}}
\newcommand{\Cov}{\mathrm{Cov}}

\title{Covering the Square Lattice and the Honeycomb\\ with Non-Overlapping Unit Disks}

\author[1,3]{Teeradej Kittipassorn}
\author[2]{Nattawut Phetmak}
\author[1]{Kirati Sriamorn}
\affil[1]{Department of Mathematics and Computer Science, Faculty of Science, Chulalongkorn University, Bangkok, Thailand}
\affil[2]{Department of Computer Engineering, Kasetsart University, Bangkok, Thailand}
\affil[3]{Centre of Excellence in Mathematics, Ministry of Higher Education, Science, Research and Innovation, Thailand}
\affil[ ]{\normalfont\small Email addresses: \texttt{teeradej.k@chula.ac.th}, \texttt{nattawut.p@ku.th}, \texttt{kirati.s@chula.ac.th}}
\date{}

\begin{document}
\maketitle

\begin{abstract}
For $d>0$, let $\Sq_d$ be the square lattice and $\Hc_d$ the honeycomb point set with minimum distance~$d$. We study for which $d$ these sets can be covered by closed unit disks with pairwise disjoint interiors; the periodic constructions of Phetmak, Nuntasombat and Fakcharoenphol cover them only from $d=1/2$ and $d=2/\sqrt{37}$ on, with gaps. Writing $\NN(d)$ for the largest number of points in a closed unit disk and $\rho$ for the number of points per unit area, we show that neither set can be covered when $\NN(d)<2\sqrt3\,\rho$, and we determine these spacings exactly, thirteen intervals for $\Sq_d$ and twelve for $\Hc_d$. At their right endpoints the bound is attained; there $\Sq_d$ is never coverable, and $\Hc_d$ is not coverable at nine of the twelve endpoints. We also give new periodic coverings, of $\Sq_d$ at $d=2/\sqrt{34}$ and of $\Hc_d$ for $2/9\le d\le1/\sqrt{19}$. Consequently, the smallest coverable spacing lies in $(0.2581,\,2/\sqrt{34}]$ for the square lattice and in $(0.1781,\,2/9]$ for the honeycomb point set, and coverability is not monotone in~$d$ for either set.
\end{abstract}



\section{Introduction}\label{sec:intro}

The \emph{coin-covering problem} asks which finite point sets in the plane can be covered by pairwise non-overlapping unit disks \cite{Win10a,Win10b,AHIU12,OKU11}. Its infinite version asks, for a periodic point set $\X$ and a spacing $d>0$, whether all points of $d\X$ can be covered by non-overlapping closed unit disks. Phetmak, Nuntasombat and Fakcharoenphol \cite{PNF26} studied this question for the triangular lattice, the square lattice and the honeycomb point set by means of periodic motif patterns. In a companion paper \cite{KPS26} we determined, for the triangular lattice, the set of spacings excluded by a density argument, decided the borderline spacings at which this argument is attained, and deduced bounds for the smallest coverable spacing. In this paper we treat the square lattice and the honeycomb point set. For both we determine exactly the spacings excluded by the density argument; at the borderline spacings the square lattice is never coverable, since no hexagonal packing is commensurable with it, while for the honeycomb point set the answer depends on an arithmetic condition. We also give new coverings, which lower the smallest known coverable spacing from $1/2$ to $2/\sqrt{34}$ for the square lattice and from $2/\sqrt{37}$ to $2/9$ for the honeycomb point set. The main tool is a general criterion for periodic point sets (Theorem~\ref{thm:criterion}). (A relaxed version of the finite problem, in which the disks may overlap but every point must be covered exactly once, was studied by Chun, Kipp and Roch \cite{CKR25}.) The smallest known finite set that cannot be covered by disjoint open unit disks, due to Nakai and Bundo \cite{NB26}, consists of a $9\times4$ piece of the square lattice $\Sq_a$ with $a=0.2531$ together with $8$ further points; they report that it becomes coverable for $a=0.2532$. For comparison, the whole lattice $\Sq_d$ is non-coverable for every $d\le0.2581$ (Theorem~\ref{thm:sq}).

Let $\Sq_d$ be the square lattice and $\Hc_d$ the honeycomb point set, both with minimum distance $d$ (Section~\ref{sec:setting}). The constructions of \cite{PNF26} give the following.

\begin{theorem}[{\cite[Theorems~4.1 and~5.2]{PNF26}}]\label{thm:PNF}
The square lattice $\Sq_d$ is coverable whenever
\[
d\in\Bigl[\tfrac12,\tfrac{6}{\sqrt{130}}\Bigr]\cup\Bigl[\tfrac{2}{\sqrt{13}},\tfrac{2}{\sqrt{10}}\Bigr]\cup\Bigl[\tfrac{4}{\sqrt{37}},\tfrac{1}{\sqrt2}\Bigr]\cup\Bigl[\tfrac{4}{\sqrt{26}},\infty\Bigr),
\]
and the honeycomb point set $\Hc_d$ is coverable whenever
\begin{gather*}
d\in\Bigl[\tfrac{2}{\sqrt{37}},\tfrac{1}{\sqrt7}\Bigr]\cup\Bigl[\tfrac{2}{\sqrt{27}},\tfrac{3}{\sqrt{57}}\Bigr]\cup\Bigl[\tfrac{2}{\sqrt{21}},\tfrac{5}{\sqrt{91}}\Bigr]\\
\cup\Bigl[\tfrac{2}{\sqrt{13}},\tfrac{3}{\sqrt{21}}\Bigr]\cup\Bigl[\tfrac{4}{\sqrt{37}},\infty\Bigr).
\end{gather*}
\end{theorem}

On the other side, the argument of Okayama, Kiyomi and Uehara \cite[Lemma~4]{OKU11} shows that $\Sq_d$ is not coverable for $d<\sqrt2\,(2/\sqrt3-1)\approx0.2188$ and $\Hc_d$ is not coverable for $d<2/\sqrt3-1\approx0.1547$. (Their proof uses only that every disk of radius $2/\sqrt3-1$ inside a certain disk contains a point of the set, which holds as soon as the covering radius, $d/\sqrt2$ for $\Sq_d$ and $d$ for $\Hc_d$, is less than $2/\sqrt3-1$; see also \cite{KPS26}.) Between these values and the constructions of Theorem~\ref{thm:PNF} the status was unknown.

\paragraph{The density obstruction.} Let $\NN(d)$ denote the largest number of points of the point set contained in a closed unit disk, and let $\rho$ be the number of points per unit area, so $\rho=1/d^2$ for $\Sq_d$ and $\rho=4/(3\sqrt3\,d^2)$ for $\Hc_d$. By the Thue--Fejes T\'oth theorem, a covering forces some disk to contain at least $2\sqrt3\,\rho$ points, so the set is not coverable whenever $\NN(d)<2\sqrt3\,\rho$ (Proposition~\ref{prop:density}). Computing $\NN$ exactly gives the following two theorems; Tables~\ref{tab:Esq} and~\ref{tab:Ehc} list the intervals with their exact endpoints.

\begin{mainthm}\label{thm:sq}
The set $\mathcal{E}_\Sq=\{d\in(0,2]:\NN(d)<2\sqrt3/d^2\}$ for the square lattice is the union of the thirteen intervals of Table~\ref{tab:Esq}. Their right endpoints are $d_k=\sqrt{2\sqrt3/k}$ for $k\in\{52,\allowbreak 44,\allowbreak 37,\allowbreak 32,\allowbreak 28,\allowbreak 27,\allowbreak 26,\allowbreak 24,\allowbreak 22,\allowbreak 21,\allowbreak 17,\allowbreak 16,\allowbreak 14\}$. The square lattice $\Sq_d$ is non-coverable for every $d\in\mathcal{E}_\Sq$ and also for every right endpoint $d=d_k$. In particular, $\Sq_d$ is non-coverable for every $d\le\sqrt{2\sqrt3/52}\approx0.2581$.
\end{mainthm}

\begin{mainthm}\label{thm:hc}
The set $\mathcal{E}_\Hc=\{d\in(0,2]:\NN(d)<8/(3d^2)\}$ for the honeycomb point set is the union of the twelve intervals of Table~\ref{tab:Ehc}. Their right endpoints are $e_k=\sqrt{8/(3k)}$ for $k\in\{84,\allowbreak 58,\allowbreak 54,\allowbreak 37,\allowbreak 34,\allowbreak 30,\allowbreak 28,\allowbreak 27,\allowbreak 26,\allowbreak 25,\allowbreak 18,\allowbreak 15\}$. The honeycomb point set $\Hc_d$ is non-coverable for every $d\in\mathcal{E}_\Hc$ and also for $d=e_k$ whenever $2k$ is not the norm of an Eisenstein integer, that is, for $k\in\{84,\allowbreak 58,\allowbreak 37,\allowbreak 34,\allowbreak 30,\allowbreak 28,\allowbreak 27,\allowbreak 25,\allowbreak 15\}$. In particular, $\Hc_d$ is non-coverable for every $d\le e_{84}=2/(3\sqrt{14})\approx0.1782$.
\end{mainthm}

Of the remaining three endpoints of $\mathcal{E}_\Hc$, the spacing $e_{18}=2/\sqrt{27}$ is coverable by Theorem~\ref{thm:PNF} and $e_{54}=2/9$ is coverable by Theorem~\ref{thm:new} below, while $e_{26}=2/\sqrt{39}$ remains open. This spacing lies between two intervals of $\mathcal{E}_\Hc$, so either $\Hc_d$ is coverable at an isolated spacing, or the non-coverable spacings near $2/\sqrt{39}$ form a single interval.

\paragraph{New coverings.} A computer search for single-family lattice patterns, in which the disk centres form a translate of a sublattice of $\Z[i]$ or $\Z[\omega]$, gives coverings that are not among those of \cite{PNF26}.

\begin{mainthm}\label{thm:new}
The square lattice $\Sq_d$ is coverable for $d\in\{2/\sqrt{34},\,1/\sqrt5,\,2/\sqrt{17},\,\sqrt{10/17}\}$ and for $d\in[1/\sqrt2,\,6/\sqrt{65}]$. The honeycomb point set $\Hc_d$ is coverable for $d\in[2/9,\,1/\sqrt{19}]$ and for $d=1/\sqrt{13}$.
\end{mainthm}

The spacings $2/\sqrt{34}$, $1/\sqrt5$ and $2/\sqrt{17}$ are left endpoints of intervals of $\mathcal E_\Sq$, and $1/\sqrt{19}$ is a left endpoint of an interval of $\mathcal E_\Hc$. At these spacings, coverability therefore changes exactly at the endpoint. The square-lattice spacings $2/\sqrt{17}$, $\sqrt{10/17}$ and the interval $[1/\sqrt2,6/\sqrt{65}]$ also follow from the inclusion $(a+bi)\Sq_d\subseteq\Sq_d$ (Lemma~\ref{lem:sub}); together with Theorem~\ref{thm:PNF} they show that $\Sq_d$ is coverable on $[4/\sqrt{37},6/\sqrt{65}]$. The covering of $\Hc_d$ at $1/\sqrt{13}$ is the restriction of a covering of the triangular lattice; see Section~\ref{sec:new}.

\paragraph{Consequences.} Since the coverable spacings form a closed set (Proposition~\ref{prop:closed}), there are smallest coverable spacings.

\begin{maincor}\label{cor:dmin}
Let $d_{\min}^\Sq$ and $d_{\min}^\Hc$ be the smallest spacings for which $\Sq_d$ and $\Hc_d$, respectively, are coverable. Then
\[
\sqrt{2\sqrt3/52}<d_{\min}^\Sq\le\frac{2}{\sqrt{34}}\qquad\text{and}\qquad\frac{2}{3\sqrt{14}}<d_{\min}^\Hc\le\frac{2}{9} .
\]
\end{maincor}

Moreover, coverability is not monotone in $d$ for either set. The honeycomb point set is coverable at $d=1/\sqrt7$ and at $d=2/\sqrt{27}$ but at no $d$ strictly between them, so the gap $(1/\sqrt7,2/\sqrt{27})$ of Theorem~\ref{thm:PNF} consists entirely of non-coverable spacings. The square lattice is coverable at $d=2/\sqrt{34}$ but not on the interval of $\mathcal E_\Sq$ that starts there, and coverable again at $d=1/\sqrt5$. All intervals of $\mathcal{E}_\Sq$ except the first lie in the gap $(0.2188,1/2)$; the three gaps of Theorem~\ref{thm:PNF} above $1/2$ are not reached by the density argument.

\begin{figure}[t]
\centering
\includegraphics[width=\linewidth]{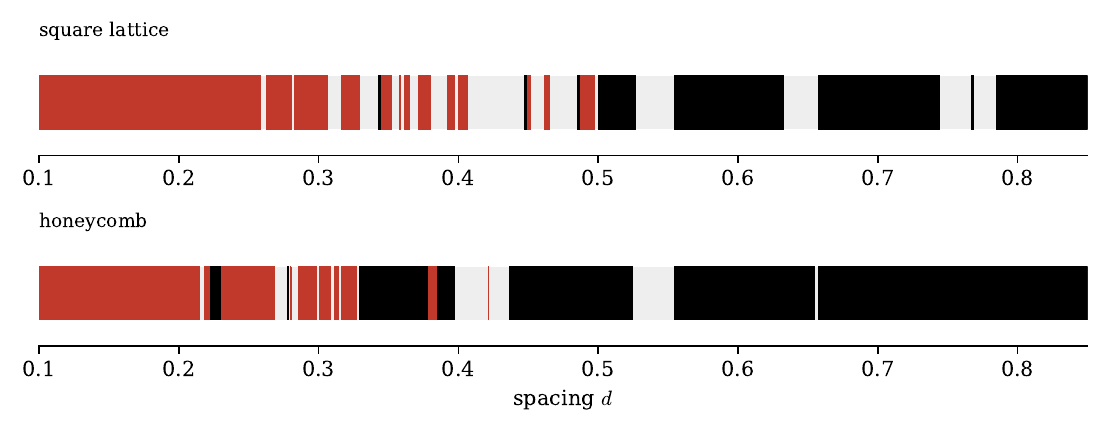}
\caption{Status of the spacings $0.1\le d\le0.85$. Black: coverable by Theorems~\ref{thm:PNF} and~\ref{thm:new} (isolated coverable spacings are drawn as thin black bars). Red: non-coverable by Theorems~\ref{thm:sq} and~\ref{thm:hc}. Grey: open.}
\label{fig:summary}
\end{figure}

\paragraph{Organisation.} Section~\ref{sec:setting} fixes notation. Section~\ref{sec:tools} develops the method for a general periodic point set, with complete proofs; the case of the triangular lattice is treated in \cite{KPS26}. Sections~\ref{sec:square} and~\ref{sec:honey} contain the non-coverability results for the two point sets, Section~\ref{sec:new} the new coverings, and Section~\ref{sec:open} lists open problems. The tables of $\NN$ are in the appendix, and the verification program is an ancillary file on the arXiv page of this paper.

\section{Setting}\label{sec:setting}

We identify $\R^2$ with $\mathbb{C}$. Let $\D(c,r)$ be the closed disk with centre $c$ and radius $r$, and $B_R=\D(0,R)$. A family of closed unit disks with centre set $C$ is \emph{non-overlapping} if $\norm{c-c'}\ge2$ for distinct $c,c'\in C$; a set is \emph{coverable} if it is contained in the union of such a family.

Let $\omega=e^{2\pi i/3}$. Following \cite{PNF26}, put
\[
\Sq=\{a+bi:a,b\in\Z\},\qquad \Hc=\{a+b\omega:a,b\in\Z,\ a+b\not\equiv0\pmod3\},
\]
and $\Sq_d=d\Sq$, $\Hc_d=d\Hc$. Both have minimum distance $d$. We write $\X$ for either set and use the following data.
\begin{center}
\renewcommand{\arraystretch}{1.25}
\begin{tabular}{@{}lll@{}}
\toprule
& square lattice $\Sq$ & honeycomb $\Hc$\\
\midrule
translation lattice $\Lambda_\X$ & $\Sq$ & $\{a+b\omega:a+b\equiv0\pmod3\}=(1-\omega)\Z[\omega]$\\
basis $v_1,v_2$ of $\Lambda_\X$ & $1,\ i$ & $1-\omega,\ 1+2\omega$\\
points per unit area $\rho_1$ & $1$ & $4/(3\sqrt3)$\\
Voronoi cell of a point & square of side $1$ & equilateral triangle of side $\sqrt3$\\
covering radius $\mu_\X$ of $\X$ & $1/\sqrt2$ & $1$\\
covering radius of $\Lambda_\X$ & $1/\sqrt2$ & $1$\\
\bottomrule
\end{tabular}
\end{center}
Here the covering radius of a set $Y$ is $\sup_{x\in\R^2}\operatorname{dist}(x,Y)$. All Voronoi cells of $\X$ are congruent, with area $1/\rho_1$. After scaling by $d$, the density is $\rho=\rho_1/d^2$ and the covering radii are multiplied by $d$. Norms are $\Nm(a+bi)=a^2+b^2$ and $\Nm(a+b\omega)=a^2-ab+b^2$. For $\alpha=a+bi$ and $\beta=x+yi$ the triangle $0,\alpha,\beta$ has area $\frac12\abs{ay-bx}$; for Eisenstein integers it has area $\frac{\sqrt3}{4}\abs{ay-bx}$.

\begin{definition}
$\NN_\X(d)=\max_{c\in\R^2}\card(\X_d\cap\D(c,1))$. We write $\NN$ when $\X$ is clear.
\end{definition}

\section{The method for a periodic point set}\label{sec:tools}

Throughout this section $\X\in\{\Sq,\Hc\}$. The same statements for the triangular lattice are proved in \cite{KPS26}. To keep the present paper self-contained, we give complete proofs in the generality needed here; they follow those of \cite{KPS26} closely, with the changes needed for a general periodic point set.

\subsection{Compactness and density}

\begin{proposition}\label{prop:closed}
The set $\Cov_\X=\{d>0:\X_d\text{ is coverable}\}$ is closed in $(0,\infty)$.
\end{proposition}

\begin{proof}
Let $d_m\in\Cov_\X$, $d_m\to d>0$, and let $C_m$ be the centre set of a covering of $\X_{d_m}$. For every $x\in\X$ choose $c_m(x)\in C_m$ with $\norm{c_m(x)-d_mx}\le1$. For fixed $x$ these points are bounded. Since $\X$ is countable, enumerate it as $x_1,x_2,\dots$; choose nested infinite index sets $I_1\supseteq I_2\supseteq\cdots$ such that $(c_m(x_j))_{m\in I_j}$ converges, and let $m_j$ be the $j$-th element of $I_j$. Along the diagonal sequence $(m_j)$ we have $c_{m_j}(x)\to c(x)$ for every $x\in\X$. Then $\norm{c(x)-dx}\le1$, and if $c(x)\ne c(x')$, the centres $c_{m_j}(x)\ne c_{m_j}(x')$ of $C_{m_j}$ are at distance at least $2$ for large $j$, so $\norm{c(x)-c(x')}\ge2$. Hence $\{c(x):x\in\X\}$ is the centre set of a covering of $\X_d$.
\end{proof}

For a covering $C$ of $\X_d$ (identified with its set of centres) let
\[
m(C)=\sup_{x\in\X_d}\operatorname{dist}(x,C)\le1,\qquad s(C)=\inf\{\norm{c-c'}:c,c'\in C,\ c\ne c'\}\ge2 .
\]
If $C$ is periodic, the supremum is attained, so $m(C)<1$ means that no point of $\X_d$ lies on the boundary of a disk.

\begin{proposition}\label{prop:slack}
Let $C$ be a covering of $\X_d$.
\begin{enumerate}[label=(\alph*),nosep]
\item If $m(C)<1$, then $\X_{d'}$ is coverable for every $d'\in[d,\,d/m(C)]$, and for $d'$ in the interior of this interval it has a covering $C'$ with $m(C')<1$ and $s(C')>2$.
\item If $s(C)>2$, then $\X_{d'}$ is coverable for every $d'\in[2d/s(C),\,d]$, and for $d'$ in the interior of this interval it has a covering $C'$ with $m(C')<1$ and $s(C')>2$.
\item If $d\in\Cov_\X$ and $(d,d+\varepsilon)\cap\Cov_\X=\emptyset$ for some $\varepsilon>0$, then every covering $C$ of $\X_d$ has $m(C)=1$, and some covering of $\X_d$ has a point of $\X_d$ that lies on the boundary of one of its disks and in the interior of none. Similarly, if $(d-\varepsilon,d)\cap\Cov_\X=\emptyset$, then every covering has $s(C)=2$.
\end{enumerate}
\end{proposition}

\begin{proof}
For $\lambda>0$ the set $C'=\lambda C$ covers $\lambda\X_d=\X_{\lambda d}$ if $\lambda m(C)\le1$, and then $m(C')=\lambda m(C)$ and $s(C')=\lambda s(C)$. Take $\lambda=d'/d$, which lies in $[1,1/m(C)]$ in case~(a) and in $[2/s(C),1]$ in case~(b); in the interior of these ranges, $\lambda m(C)<1$ and $\lambda s(C)>2$.

(c) The first and last statements follow from (a) and (b). For the second, let $C$ be a covering and $x_n\in\X_d$ with $\operatorname{dist}(x_n,C)\to1$. Choose $\ell_n\in d\Lambda_\X$ such that $x_n-\ell_n\in dR_\X$; passing to a subsequence, $x_n-\ell_n=p_0$ is constant. Then $C_n=C-\ell_n$ is a covering of $\X_d$ with $\operatorname{dist}(p_0,C_n)\to1$. For $p\in\X_d$ let $c_n(p)$ be a point of $C_n$ nearest to $p$. The diagonal argument in the proof of Proposition~\ref{prop:closed} gives a subsequence with $c_n(p)\to c(p)$ for all $p$, and $C^\ast=\{c(p)\}$ is a covering of $\X_d$. Then $\norm{c(p_0)-p_0}=1$, and every $c(p)\in C^\ast$ satisfies $\norm{c(p)-p_0}=\lim\norm{c_n(p)-p_0}\ge1$. So $p_0$ lies on the boundary of a disk of $C^\ast$ and in the interior of none.
\end{proof}

\begin{lemma}\label{lem:count}
For every $x\in\R^2$ and $R>d\mu_\X$, $\card(\X_d\cap\D(x,R))\ge\rho\,\pi(R-d\mu_\X)^2$.
\end{lemma}

\begin{proof}
The Voronoi cells of the points of $\X_d$ have area $1/\rho$, and every point of $\D(x,R-d\mu_\X)$ lies in the cell of a point within distance $d\mu_\X$, hence in $\D(x,R)$. So the cells of the points of $\X_d\cap\D(x,R)$ cover $\D(x,R-d\mu_\X)$.
\end{proof}

\begin{proposition}\label{prop:density}
If $\X_d$ is coverable, then $\NN(d)\ge2\sqrt3\,\rho=2\sqrt3\,\rho_1/d^2$.
\end{proposition}

\begin{proof}
Let $C$ be the centre set of a covering. Every point of $\X_d\cap B_R$ lies in a disk with centre in $B_{R+1}$, and each disk contains at most $\NN(d)$ points, so $\card(\X_d\cap B_R)\le\NN(d)\,\card(C\cap B_{R+1})$. Divide by $\pi R^2$ and let $R\to\infty$. By Lemma~\ref{lem:count} the left-hand side tends to $\rho$, and by the theorem of Thue and Fejes T\'oth \cite{Thu10,FT40} a set of points with pairwise distances at least $2$ satisfies $\limsup_{R\to\infty}\card(C\cap B_R)/(\pi R^2)\le1/(2\sqrt3)$. Hence $\rho\le\NN(d)/(2\sqrt3)$.
\end{proof}

\subsection{Structure of \texorpdfstring{$\NN$}{N}}

\begin{lemma}\label{lem:two}
Let $S\subset\R^2$ be a finite set with at least two points, and suppose $S\subseteq\D(c_0,r)$ for some $c_0$. Then there are distinct $p,q\in S$ and a point $c$ with $\norm{c-p}=\norm{c-q}=r$ and $S\subseteq\D(c,r)$.
\end{lemma}

\begin{proof}
The set $\bigcap_{s\in S}\D(s,r)$ of centres of radius-$r$ disks containing $S$ is nonempty and compact, so it has a boundary point $c_1$, at which some constraint is active: $\norm{c_1-p}=r$ for some $p\in S$. Parametrise the circle $\partial\D(p,r)$ by $c(\theta)=p+r(\cos\theta,\sin\theta)$ and let $\Theta=\{\theta:\norm{c(\theta)-s}\le r\text{ for all }s\in S\}$. This set is closed and contains the angle of $c_1$. It is not the whole circle: for $q\in S\setminus\{p\}$, the point $p+r\frac{p-q}{\norm{p-q}}$ is at distance $r+\norm{p-q}>r$ from $q$. Hence $\Theta$ has a boundary point $\theta^\ast\in\Theta$. If all constraints with $s\ne p$ were strict at $c(\theta^\ast)$, they would remain strict near $\theta^\ast$, contradicting $\theta^\ast\in\partial\Theta$. So $\norm{c(\theta^\ast)-q}=r$ for some $q\in S\setminus\{p\}$, and $c=c(\theta^\ast)$ has the required properties.
\end{proof}

Let $R_\X$ be a set of representatives of $\X$ modulo $\Lambda_\X$: $R_\Sq=\{0\}$ and $R_\Hc=\{1,-1\}$.

\begin{corollary}\label{cor:formula}
Let $r=1/d\ge\frac12$. Then $\NN(d)$ is the maximum of $\card(\X\cap\D(c,r))$ over all $p\in R_\X$, $q\in\X$ with $0<\abs{q-p}\le2r$ and the (one or two) centres $c$ with $\abs{c-p}=\abs{c-q}=r$.
\end{corollary}

\begin{proof}
Since the minimum distance of $\X$ is $1\le2r$, an optimal disk of radius $r$ contains at least two points of $\X$. By Lemma~\ref{lem:two} we may assume that two points $p',q'$ of $\X$ lie on its boundary. A translation by an element of $\Lambda_\X$ maps $\X$ to itself and $p'$ to a representative $p\in R_\X$.
\end{proof}

\begin{lemma}\label{lem:monotone}
$\NN$ is non-increasing and left-continuous on $(0,\infty)$.
\end{lemma}

\begin{proof}
Let $g(r)$ be the largest number of points of $\X$ in a closed disk of radius $r$, so that $\NN(d)=g(1/d)$. It suffices to show that $g$ is non-decreasing, which is clear, and right-continuous. Fix $r$; since $g$ is integer-valued and non-decreasing, there are $k$ and $\delta>0$ with $g(r')=k$ for $r'\in(r,r+\delta)$. Take $r_j\downarrow r$ and disks $\D(c_j,r_j)$ containing sets $S_j\subset\X$ with $\card S_j=k$. After a translation by an element of $\Lambda_\X$ we may assume that $S_j$ contains a point of the finite set $R_\X$, so $S_j$ lies in a fixed bounded set. Passing to a subsequence, $S_j=S$ is constant and $c_j\to c$. Then $\norm{s-c}\le r$ for all $s\in S$, so $g(r)\ge k$.
\end{proof}

The discontinuities of $\NN$ are governed by smallest enclosing disks; see \cite{Wel91}.

\begin{lemma}\label{lem:sed}
Let $S\subset\R^2$ be finite with $\card S\ge2$, and let $r(S)$ be the smallest radius of a closed disk containing $S$. Then either $r(S)=\frac12\norm{p-q}$ for some $p,q\in S$, or $r(S)$ is the circumradius of an acute triangle with vertices in $S$.
\end{lemma}

\begin{proof}
A smallest enclosing disk $\D(c,r_0)$ exists by compactness. Let $B=S\cap\partial\D(c,r_0)$. Then $c\in\operatorname{conv}B$: otherwise there is a unit vector $u$ with $\langle u,b-c\rangle>0$ for all $b\in B$, and for small $\varepsilon>0$ the disk centred at $c+\varepsilon u$ contains $S$ in its interior, so a smaller disk contains $S$. By Carath\'eodory's theorem, $c$ lies on a segment $[p,q]$ with $p,q\in B$, in which case $c=\frac{p+q}2$ and $r_0=\frac12\norm{p-q}$, or in the interior of a triangle with vertices in $B$, in which case $c$ is its circumcentre and the triangle is acute.
\end{proof}

\begin{lemma}\label{lem:jumps}
If $\NN$ is discontinuous at $d^\ast$, then $(d^\ast)^2$ equals $4/\Nm(\alpha)$ for some nonzero difference $\alpha$ of two points of $\X$, or
\[
(d^\ast)^2=\frac{4\,(ay-bx)^2}{\Nm(\alpha)\Nm(\beta)\Nm(\beta-\alpha)}\ \ (\X=\Sq),\qquad (d^\ast)^2=\frac{3\,(ay-bx)^2}{\Nm(\alpha)\Nm(\beta)\Nm(\beta-\alpha)}\ \ (\X=\Hc),
\]
where $\alpha=a+bi$, $\beta=x+yi$ (respectively $a+b\omega$, $x+y\omega$) are differences $q-p$, $q'-p$ of points $p,q,q'$ of $\X$ with $ay-bx\ne0$. If $d^\ast\ge d_0$, all norms involved are at most $4/d_0^2$, so the possible discontinuities in $[d_0,\infty)$ are finitely many.
\end{lemma}

\begin{proof}
By Lemma~\ref{lem:monotone}, a discontinuity at $d^\ast$ means $\NN(d)<\NN(d^\ast)=:k$ for $d>d^\ast$. Then $k\ge2$, and some $k$-point set $S\subset\X$ lies in a disk of radius $r^\ast=1/d^\ast$ but in no smaller disk, i.e.\ $r(S)=r^\ast$. By Lemma~\ref{lem:sed}, $r^\ast$ is half the distance $\sqrt{\Nm(\alpha)}$ of two points of $S$, or the circumradius of a triangle with vertices $p,q,q'\in S$. In the second case the formulas follow from $R=\frac{\text{product of sides}}{4\cdot\text{area}}$ and the area formulas of Section~\ref{sec:setting}, with $\alpha=q-p$, $\beta=q'-p$. All points of $S$ lie in a disk of radius $r^\ast\le1/d_0$, so all the norms are at most $4/d_0^2$.
\end{proof}

Consequently, $\NN$ is constant on each interval between consecutive possible discontinuities, and its value there can be computed at one interior point by Corollary~\ref{cor:formula}.

\begin{lemma}\label{lem:crude}
Let $\delta_\X=\bigl(\sqrt{2\sqrt3/\pi}-1\bigr)/\mu_\X$, that is, $\delta_\Sq\approx0.07082$ and $\delta_\Hc\approx0.05008$. If $0<d\le\delta_\X$, then $\NN(d)<2\sqrt3\rho_1/d^2$.
\end{lemma}

\begin{proof}
If $\D(c,1)$ contains $n$ points of $\X_d$, their Voronoi cells have disjoint interiors, area $d^2/\rho_1$ and circumradius $d\mu_\X$, so they lie in $\D(c,1+d\mu_\X)$. Hence $n\le\rho_1\pi(1+d\mu_\X)^2/d^2$, and $\pi(1+d\mu_\X)^2\le2\sqrt3$ for $d\le\delta_\X$, with strict inequality for the count since finitely many cells cannot fill a disk.
\end{proof}

\begin{lemma}[Chain]\label{lem:chain}
If $\delta_0<\dots<\delta_M$ and $n_i\ge\NN(\delta_i)$ are integers with $\delta_{i+1}=\sqrt{2\sqrt3\rho_1/n_i}$, then $\NN(d)<2\sqrt3\rho_1/d^2$ for all $d\in[\delta_0,\delta_M)$.
\end{lemma}

\begin{proof}
For $d\in[\delta_i,\delta_{i+1})$, monotonicity gives $\NN(d)\le\NN(\delta_i)\le n_i=2\sqrt3\rho_1/\delta_{i+1}^2<2\sqrt3\rho_1/d^2$.
\end{proof}

\subsection{Borderline spacings}

Let $C\subset\R^2$ be a set with $\norm{c-c'}\ge2$ for distinct $c,c'\in C$, and let $V_c=\{x:\norm{x-c}\le\norm{x-c'}\text{ for all }c'\in C\}$ be the Voronoi cell of $c\in C$. Since $\norm{x-c'}\ge\norm{c-c'}-\norm{x-c}\ge1$ whenever $\norm{x-c}\le1$, we have $\D(c,1)\subseteq V_c$.

\begin{lemma}\label{lem:cell}
If $V_c$ is bounded, then $\abs{V_c}\ge2\sqrt3$, with equality if and only if $V_c$ is a regular hexagon with inradius $1$ centred at $c$.
\end{lemma}

\begin{proof}
Translate so that $c=0$. The bounded cell $V_c$ is a convex polygon. Each edge $e$ lies on the perpendicular bisector of $0$ and some $c'\in C$, at distance $h_e=\norm{c'}/2\ge1$ from $0$. Each vertex $v$ is equidistant from $0$ and at least two further points of $C$; these three points lie on a circle of radius $\norm v$ and are pairwise at distance at least $2$, so the triangle they form has an angle $A\le\pi/3$ and $2\le2\norm v\sin A\le\sqrt3\,\norm v$, i.e.\ $\norm v\ge2/\sqrt3$.

Let $T_e=\operatorname{conv}(\{0\}\cup e)$ and let $\varphi_e$ be its angle at $0$, so that $\sum_e\varphi_e=2\pi$. If $u,w$ are the endpoints of $e$, then $\abs{T_e}=\frac12\norm u\norm w\sin\varphi_e\ge\frac23\sin\varphi_e$. Also, with $a,b$ the signed angles from the foot of the perpendicular from $0$ to $u$ and $w$, we have $\abs{T_e}=\frac{h_e^2}{2}\abs{\tan b-\tan a}$ with $\abs{b-a}=\varphi_e$; if $a\le0\le b$, convexity of $\tan$ on $[0,\pi/2)$ gives $\abs{\tan b-\tan a}\ge2\tan(\varphi_e/2)$, and if $0\le a<b$, then $\tan b-\tan a=\frac{\sin\varphi_e}{\cos a\cos b}\ge\frac{2\sin\varphi_e}{1+\cos\varphi_e}=2\tan(\varphi_e/2)$. Hence $\abs{T_e}\ge\gamma(\varphi_e)$ with $\gamma(\varphi)=\max\{\frac23\sin\varphi,\tan(\varphi/2)\}$. Since $\sin\varphi/\varphi$ is decreasing and $\tan(\varphi/2)/\varphi$ is increasing, $\gamma(\varphi)\ge\frac{\sqrt3}{\pi}\varphi$, with equality only at $\varphi=\pi/3$. Summing over the edges gives $\abs{V_c}\ge2\sqrt3$. Equality forces $\varphi_e=\pi/3$ and $\norm u=\norm w=2/\sqrt3$ for every edge, i.e.\ a regular hexagon with inradius $1$.
\end{proof}

Suppose now that $\NN(d)\le k=2\sqrt3\,\rho$, and let ${R_0}=1+d\mu_\X$. For a covering $C$ of $\X_d$ and $c\in C$ put
\[
n_c=\card(\X_d\cap\D(c,1)),\qquad \delta(c)=\rho\abs{V_c}-n_c .
\]

\begin{lemma}\label{lem:defect}
Let $C$ be a covering of $\X_d$. Then $V_c\subseteq\D(c,{R_0})$ and $\delta(c)\ge0$ for every $c\in C$, and there is a constant $K$ such that $\sum_{c\in C\cap\D(x,R)}\delta(c)\le KR$ for all $x\in\R^2$ and $R\ge1$.
\end{lemma}

\begin{proof}
Every point of the plane is within $d\mu_\X$ of a point of $\X_d$, which is within $1$ of a centre; so $x\in V_c$ implies $\norm{x-c}\le{R_0}$. By Lemma~\ref{lem:cell} and $n_c\le\NN(d)\le k$, $\delta(c)\ge2\sqrt3\rho-k=0$. The cells of the centres in $\D(x,R)$ lie in $\D(x,R+{R_0})$ and have disjoint interiors, so $\sum\abs{V_c}\le\pi(R+{R_0})^2$. Every point of $\X_d\cap\D(x,R-1)$ lies in a disk with centre in $\D(x,R)$, so $\sum n_c\ge\rho\pi(R-1-d\mu_\X)^2$ by Lemma~\ref{lem:count}. Subtracting gives $\sum\delta(c)\le KR$ for a suitable $K$.
\end{proof}

\begin{lemma}\label{lem:perfect}
If $\NN(d)=k=2\sqrt3\,\rho$ and $\X_d$ is coverable, then there is a covering $C^\ast$ of $\X_d$ whose centre set is $c_0+L'$ for a triangular lattice $L'$ of minimum distance $2$, and in which every disk contains exactly $k$ points of $\X_d$.
\end{lemma}

\begin{proof}
Let $C$ be a covering. Discarding disks without points of $\X_d$, we may assume every centre covers a point. Two centres at distance at least $2$ in $\D(p,1)$ are antipodal, so every point of $\X_d$ lies in at most two disks.

\emph{Step 1: almost perfect patches.} Let $R\ge16$ and $r=R^{1/4}$. The disks $\D(y,r)$ with $y\in2r\Z^2\cap\D(0,R-r)$ lie in $\D(0,R)$, every point lies in at most two of them, and there are at least $\beta R^2/r^2$ of them for an absolute constant $\beta>0$ and large $R$. By Lemma~\ref{lem:defect} one of them satisfies $\sum_{c\in C\cap\D(y,r)}\delta(c)\le\frac{2KRr^2}{\beta R^2}=\frac{2K}{\beta}R^{-1/2}$. Taking $R=R_m\to\infty$ we obtain $y_m$, $r_m\to\infty$ and $\varepsilon_m\to0$ with $\sum_{c\in C\cap\D(y_m,r_m)}\delta(c)\le\varepsilon_m$. Choose $\ell_m\in d\Lambda_\X$ with $\norm{y_m-\ell_m}$ at most the covering radius of $d\Lambda_\X$, which equals $d\mu_\X\le1$ for both point sets (Section~\ref{sec:setting}). Since $\X_d-\ell_m=\X_d$, the set $C_m=C-\ell_m$ is again a covering, and its defects satisfy $\sum_{c\in C_m\cap\D(0,r_m-1)}\delta_m(c)\le\varepsilon_m$.

\emph{Step 2: a limit covering.} For each $m$ and $p\in\X_d$, let $c^1_m(p),c^2_m(p)$ be the (one or two) centres of $C_m$ whose disks contain $p$, repeated if there is only one. By the diagonal argument in the proof of Proposition~\ref{prop:closed}, a subsequence satisfies $c^i_m(p)\to c^i_\ast(p)$ for all $p,i$. Let $C^\ast$ be the set of limits. As in that proof, $C^\ast$ is a covering of $\X_d$. If $c^i_\ast(p)=c^j_\ast(q)$, then $c^i_m(p)$ and $c^j_m(q)$ are equal or at distance at least $2$ and have the same limit, so they are equal for large $m$; if $c^i_\ast(p)\ne c^j_\ast(q)$, they are distinct for large $m$.

\emph{Step 3: semicontinuity of the defect.} Fix $c=c^i_\ast(p_0)\in C^\ast$ and let $F$ be the set of centres of $C^\ast$ within distance $2{R_0}+2$ of $c$. Then $V_c$ equals the cell $V^F_c$ of $c$ with respect to $F$: clearly $V_c\subseteq V^F_c$; if $V^F_c\not\subseteq\D(c,{R_0})$, convexity gives $x\in V^F_c$ with ${R_0}<\norm{x-c}<{R_0}+1$, and since $x\notin V_c$ by Lemma~\ref{lem:defect}, some $c''\notin F$ has $\norm{x-c''}<{R_0}+1$, contradicting $\norm{x-c''}>(2{R_0}+2)-({R_0}+1)$; and every $x\in\D(c,{R_0})$ is closer to $c$ than to any $c''\notin F$. Every centre in $F$ covers a point of the finite set $P_1=\X_d\cap\D(p_0,2{R_0}+4)$, and so does every centre of $C_m$ near $c_m=c^i_m(p_0)$. By Step~2, the finite families $\{c^j_m(q):q\in P_1\}$ converge to $\{c^j_\ast(q):q\in P_1\}$ with the same coincidences, so $\abs{V_{c_m}}\to\abs{V_c}$. Since the disks are closed, $n_c\ge\limsup_m n_{c_m}$. Hence $\delta(c)\le\liminf_m\delta_m(c_m)\le\lim\varepsilon_m=0$, using $c_m\in\D(0,r_m-1)$ for large $m$.

\emph{Step 4: rigidity.} For $c\in C^\ast$ we get $0=\delta(c)\ge\rho\abs{V_c}-k\ge0$, so $\abs{V_c}=2\sqrt3$ and $n_c=k$. By Lemma~\ref{lem:cell} every cell is a regular hexagon of inradius $1$. Adjacent cells share a full edge and are mirror images in it, so all cells have the same orientation, the neighbours of $c$ are $c+2u$ for the six unit normals $u$ of its edges, and the cells of $c_0+L'$ tile the plane, where $L'$ is the triangular lattice of minimum distance $2$ in this orientation. A further centre would have a cell overlapping these, so $C^\ast=c_0+L'$.
\end{proof}

Let $\Omega_k=\{x\in\R^2:\card(\X_d\cap\D(x,1))\ge k\}$; this set is closed, since the count is upper semicontinuous in $x$, and it is invariant under $d\Lambda_\X$.

\begin{lemma}\label{lem:line}
Let $L'$ be a triangular lattice of minimum distance $2$ such that $L'+d\Lambda_\X$ is not discrete. If $c_0+L'\subseteq\Omega_k$, then $\Omega_k$ contains a straight line.
\end{lemma}

\begin{proof}
The closure $G$ of $L'+d\Lambda_\X$ is a non-discrete closed subgroup of $\R^2$, so it contains a line through $0$ \cite[Ch.~VII, \S1, Thm.~2]{Bou98}, and $c_0+G\subseteq\Omega_k$ because $\Omega_k$ is closed and $d\Lambda_\X$-invariant.
\end{proof}

\begin{lemma}\label{lem:comm}
Let $L'$ be a triangular lattice of minimum distance $2$.
\begin{enumerate}[label=(\alph*),nosep]
\item $L'+d\Lambda_\Sq$ is never discrete.
\item If $d=e_k=\sqrt{8/(3k)}$ and $2k$ is not the norm of an Eisenstein integer, then $L'+d\Lambda_\Hc$ is not discrete.
\end{enumerate}
\end{lemma}

\begin{proof}
If $L'+d\Lambda_\X$ is discrete, it is a lattice containing $d\Lambda_\X$ with some finite index $m$, so $mL'\subseteq d\Lambda_\X$.

(a) Then $mL'\subseteq d\Z^2$ contains an equilateral triangle with vertices in $d\Z^2$. Its area would be a rational multiple of $d^2$ by the shoelace formula, and also $\frac{\sqrt3}{4}s^2$ with $s^2$ a rational multiple of $d^2$, which is impossible since $\sqrt3$ is irrational.

(b) Take $u\in L'$ with $\norm u=2$. Then $mu=d(1-\omega)\beta$ for some $\beta\in\Z[\omega]$, so $4m^2=d^2\cdot3\Nm(\beta)=8\Nm(\beta)/k$, i.e.\ $2km^2=4\Nm(\beta)=\Nm(2\beta)$. Since a positive integer is a norm if and only if every prime $\equiv2\pmod3$ divides it to an even power \cite{IR90}, $2km^2$ is a norm only if $2k$ is.
\end{proof}

\begin{lemma}[Fence]\label{lem:fence}
If there is $s\in\R^2$ such that $\card(\X_d\cap\D(x,1))\le k-1$ for all $x$ on the segments $[s,s+dv_1]$ and $[s,s+dv_2]$, then $\Omega_k$ contains no straight line.
\end{lemma}

\begin{proof}
The parallelograms $s+\lambda+\{a\,dv_1+b\,dv_2:0\le a,b\le1\}$, $\lambda\in d\Lambda_\X$, tile the plane, and their edges are $d\Lambda_\X$-translates of the two segments. A line leaves every bounded parallelogram it meets, so it crosses an edge, and by $d\Lambda_\X$-invariance of the count that point of the line is not in $\Omega_k$.
\end{proof}

\begin{theorem}\label{thm:criterion}
Let $\NN(d)=k=2\sqrt3\,\rho$. Suppose that $L'+d\Lambda_\X$ is not discrete for every triangular lattice $L'$ of minimum distance $2$, and that a fence as in Lemma~\ref{lem:fence} exists. Then $\X_d$ is not coverable.
\end{theorem}

\begin{proof}
If $\X_d$ were coverable, Lemma~\ref{lem:perfect} would give $c_0+L'\subseteq\Omega_k$, and Lemma~\ref{lem:line} a line in $\Omega_k$, contradicting Lemma~\ref{lem:fence}.
\end{proof}

\begin{remark}[Exact arithmetic]\label{rem:float}
All computations are carried out by the program \texttt{verify\_square\_\allowbreak honeycomb.py}, which is available as an ancillary file on the arXiv page of this paper; it requires Python~3 and \texttt{numpy}, used only for integer arrays, and runs in about one minute. It computes the blocks of Tables~\ref{tab:blocksq} and~\ref{tab:blockhc} (Part~A), the chains (Part~B), the random test (Part~C), the fences of Table~\ref{tab:fence} (Part~D), the coverings of Table~\ref{tab:new} (Part~E) and the sublattice transfer of Section~\ref{sec:new} (Part~F). All computations are exact: no decision uses floating-point arithmetic, and there are no tolerances. Candidate discontinuities (Lemma~\ref{lem:jumps}) are rational values of $d^2$. Every squared spacing that occurs lies in $\mathbb{Q}(\sqrt3)$, and it is rational except for the square-lattice spacings with $d^2=2\sqrt3/n$, namely the borderline spacings $d_k$ and the spacings $\delta_{i+1}$ of the chain; every decision is the sign of an integer, of an element of $\mathbb{Q}(\sqrt3)$, or of a number $x+y\sqrt z$ with $x,y,z\in\mathbb{Q}(\sqrt3)$ and $z\ge0$. The sign of $x+y\sqrt z$ is the common sign of $x$ and $y$ if these do not have opposite signs, and $\operatorname{sign}(x)\operatorname{sign}(x^2-y^2z)$ otherwise, and the sign of $a+b\sqrt3$ with $a,b\in\mathbb{Q}$ is decided in the same way.

\emph{Values of $\NN$.} These are only computed at spacings with $d^2=r^{-2}\in\mathbb{Q}$. Between two consecutive candidate discontinuities we use the rational value of $d^2$ with the smallest denominator, and the random test uses $1000$ random rational values of $d^2$ for each point set. In Corollary~\ref{cor:formula}, write $\alpha=q-p=a_1+b_1\tau$ and $\beta=x-p=a_0+b_0\tau$ for a point $x\in\X$, where $\tau=i$ or $\omega$, and let $c_\pm$ be the two centres. Since $\abs{c_\pm-p}=r$,
\[
2\bigl(\abs{x-c_\pm}^2-r^2\bigr)=\bigl(2\Nm(\beta)-2\langle\alpha,\beta\rangle\bigr)\mp\Delta\sqrt{\kappa\Bigl(\frac{4r^2}{\Nm(\alpha)}-1\Bigr)},\qquad\Delta=a_1b_0-b_1a_0\in\Z,
\]
with $\kappa=4$ for $\Sq$ and $\kappa=3$ for $\Hc$. The first term is an integer, so $x\in\D(c_\pm,r)$ is decided in exact integer arithmetic; the program evaluates these signs for all $x$ at once on integer arrays, using 64-bit integers only after checking that no overflow can occur and Python integers otherwise. It suffices to let $\alpha$ range over a sector of angle $\pi/4$ for $\Sq$ and $\pi/3$ for $\Hc$: the set $\X-p$ is invariant under the rotations by $\pi/2$ and reflections in the coordinate axes for $\Sq$, and under multiplication by $\omega$ and complex conjugation for $\Hc$, since both maps fix the residue of $a+b\omega$ modulo $1-\omega$, that is, the value of $a+b\bmod3$.

\emph{Chains.} For the honeycomb all spacings $\delta_i$ have rational squares. For the square lattice $\delta_{i+1}^2=2\sqrt3/n_i$ is irrational, and the program takes $n_{i+1}=\NN(\delta_{i+1}')$ for a rational $\delta_{i+1}'^2\le\delta_{i+1}^2$ (obtained from a convergent $p/q<\sqrt3$ with $q\le10^6$). By Lemma~\ref{lem:monotone}, $n_{i+1}\ge\NN(\delta_{i+1})$, which is all that Lemma~\ref{lem:chain} requires. The inequality $\pi(1+\mu_\X\delta_0)^2<2\sqrt3$ of Lemma~\ref{lem:crude} is checked using $\pi<355/113$.

\emph{Borderline spacings and fences.} The equality $\NN(d_k)=k$ (respectively $\NN(e_k)=k$) follows from Lemmas~\ref{lem:monotone} and~\ref{lem:jumps} once $d_k^2$ is shown to lie in a block with value $k$, which is an exact comparison in $\mathbb{Q}(\sqrt3)$; for the honeycomb the program also evaluates $\NN(e_k)$ directly. On a segment $s+\tau v$, $0\le\tau\le1$, with $s,v$ rational in the basis of $\Lambda_\X$, the set of $\tau$ for which a given point lies in the closed unit disk is an interval with endpoints $x\pm\sqrt z$, $x\in\mathbb{Q}$, $z\in\mathbb{Q}(\sqrt3)$. A sweep sorts these endpoints by exact comparisons, which reduce to signs of numbers $x+y\sqrt z$ as above, with starts before ends at equal parameters, and returns the exact maximum count over closed disks centred on the segment.

An earlier floating-point implementation of the same program, with a relative tolerance of $10^{-9}$, gives identical output.
\end{remark}

\section{The square lattice}\label{sec:square}

Here $\rho=1/d^2$, the density threshold is $2\sqrt3/d^2$, and $d_k=\sqrt{2\sqrt3/k}$.

\begin{proof}[Proof of Theorem~\ref{thm:sq}]
By Lemma~\ref{lem:crude}, $\NN(d)<2\sqrt3/d^2$ for $d\le0.0708<\delta_\Sq$. The chain of Lemma~\ref{lem:chain} with $30$ steps, from $0.0708$ to $0.20188$, covers $[0.0708,0.2)$ (Part~B of the program). On $[1/5,2]$, Lemma~\ref{lem:jumps} gives $470$ possible discontinuities, and the evaluation of $\NN$ between them gives the $51$ blocks of Table~\ref{tab:blocksq}. On a block with value $k$, $\NN(d)<2\sqrt3/d^2$ if and only if $d<d_k$, which yields Table~\ref{tab:Esq}. Non-coverability on $\mathcal{E}_\Sq$ follows from Proposition~\ref{prop:density}.

At each right endpoint $d_k$ we have $\NN(d_k)=k=2\sqrt3\rho$. By Lemma~\ref{lem:comm}(a), $L'+d\Lambda_\Sq$ is never discrete, and Part~D of the program finds a fence for every $k$ (Table~\ref{tab:fence}). Theorem~\ref{thm:criterion} shows that $\Sq_{d_k}$ is not coverable.
\end{proof}

\begin{table}[t]
\centering
\caption{The thirteen intervals of $\mathcal{E}_\Sq$, with the value $k$ of $\NN$ just left of the right endpoint $d_k=\sqrt{2\sqrt3/k}$. By Theorem~\ref{thm:sq} the right endpoints are also non-coverable.}
\label{tab:Esq}
\renewcommand{\arraystretch}{1.3}
\begin{tabular}{@{}rllr@{}}
\toprule
& interval of $\mathcal{E}_\Sq$ & decimal & $k$\\
\midrule
1 & $(0,\,\sqrt{2\sqrt3/52})$ & $(0.000000,\,0.258103)$ & 52 \\
2 & $(\tfrac{2}{\sqrt{58}},\,\sqrt{2\sqrt3/44})$ & $(0.262613,\,0.280588)$ & 44 \\
3 & $(\tfrac{2}{5\sqrt{2}},\,\sqrt{2\sqrt3/37})$ & $(0.282843,\,0.305981)$ & 37 \\
4 & $(\tfrac{1}{\sqrt{10}},\,\sqrt{2\sqrt3/32})$ & $(0.316228,\,0.329019)$ & 32 \\
5 & $(\tfrac{2}{\sqrt{34}},\,\sqrt{2\sqrt3/28})$ & $(0.342997,\,0.351736)$ & 28 \\
6 & $(\tfrac{4}{5\sqrt{5}},\,\sqrt{2\sqrt3/27})$ & $(0.357771,\,0.358190)$ & 27 \\
7 & $(\tfrac{38}{5\sqrt{442}},\,\sqrt{2\sqrt3/26})$ & $(0.361495,\,0.365013)$ & 26 \\
8 & $(\tfrac{2}{\sqrt{29}},\,\sqrt{2\sqrt3/24})$ & $(0.371391,\,0.379918)$ & 24 \\
9 & $(\tfrac{2}{\sqrt{26}},\,\sqrt{2\sqrt3/22})$ & $(0.392232,\,0.396811)$ & 22 \\
10 & $(\tfrac{2}{5},\,\sqrt{2\sqrt3/21})$ & $(0.400000,\,0.406149)$ & 21 \\
11 & $(\tfrac{1}{\sqrt{5}},\,\sqrt{2\sqrt3/17})$ & $(0.447214,\,0.451410)$ & 17 \\
12 & $(\tfrac{6}{13},\,\sqrt{2\sqrt3/16})$ & $(0.461538,\,0.465302)$ & 16 \\
13 & $(\tfrac{2}{\sqrt{17}},\,\sqrt{2\sqrt3/14})$ & $(0.485071,\,0.497429)$ & 14 \\
\bottomrule
\end{tabular}
\end{table}

\begin{proof}[Proof of Corollary~\ref{cor:dmin}, square lattice]
$\Cov_\Sq$ is closed (Proposition~\ref{prop:closed}), nonempty and contained in $(d_{52},\infty)$ by Theorem~\ref{thm:sq}. Its minimum exists, it exceeds $d_{52}$, and $2/\sqrt{34}\in\Cov_\Sq$ by Theorem~\ref{thm:new}.
\end{proof}

\begin{figure}[t]
\centering
\includegraphics[width=0.92\linewidth]{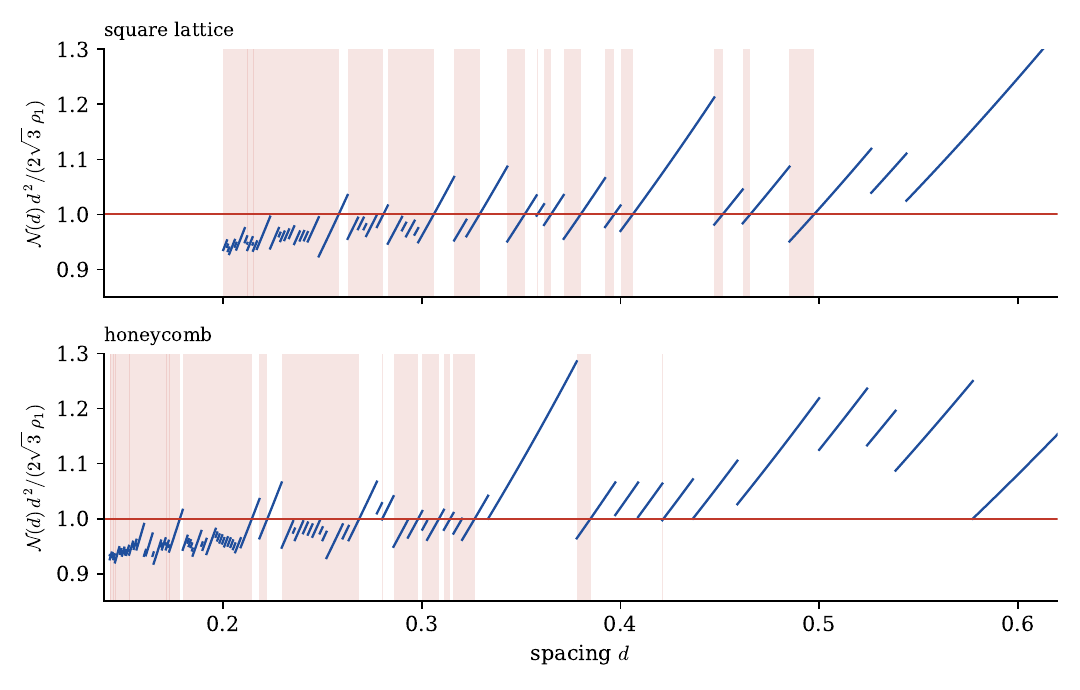}
\caption{The ratio $\NN(d)\,d^2/(2\sqrt3\rho_1)$ for the square lattice (from $d=0.2$) and the honeycomb point set (from $d=1/7$). Shaded spacings, where the ratio is below $1$, form $\mathcal E_\Sq$ and $\mathcal E_\Hc$.}
\label{fig:ratio}
\end{figure}

The last interval of $\mathcal{E}_\Sq$ ends at $d_{14}\approx0.4974$, just below the zipper construction of \cite{PNF26}, which covers $\Sq_d$ for $1/2\le d\le6/\sqrt{130}$ with $14$ points per disk. Between consecutive intervals of $\mathcal{E}_\Sq$ we have $\NN(d)\ge2\sqrt3/d^2$, and the density argument gives no information there. On the three intervals above $1/2$ that are still open, $(6/\sqrt{130},2/\sqrt{13})$, $(2/\sqrt{10},4/\sqrt{37})$ and $(6/\sqrt{65},4/\sqrt{26})$, the infimum of the ratio $\NN(d)d^2/(2\sqrt3)$ is approximately $1.025$, $1.039$ and $1.119$ (Figure~\ref{fig:ratio}).

\section{The honeycomb point set}\label{sec:honey}

Here $\rho=4/(3\sqrt3d^2)$, the density threshold is $8/(3d^2)$, and $e_k=\sqrt{8/(3k)}$.

\begin{proof}[Proof of Theorem~\ref{thm:hc}]
By Lemma~\ref{lem:crude}, $\NN(d)<8/(3d^2)$ for $d\le0.05<\delta_\Hc$. A chain of $29$ steps covers $[0.05,1/7)$. On $[1/7,2]$, Lemma~\ref{lem:jumps} gives $1388$ possible discontinuities and $85$ blocks (Table~\ref{tab:blockhc}), which yield Table~\ref{tab:Ehc}. Non-coverability on $\mathcal{E}_\Hc$ follows from Proposition~\ref{prop:density}.

At the right endpoints $e_k$ we have $\NN(e_k)=k=8/(3e_k^2)$. For $k\in\{84,\allowbreak 58,\allowbreak 37,\allowbreak 34,\allowbreak 30,\allowbreak 28,\allowbreak 27,\allowbreak 25,\allowbreak 15\}$, the number $2k$ is divisible to an odd power by a prime $\equiv2\pmod3$ (namely $2,29,2,17,5,2,2,2,2$), so it is not a norm; Lemma~\ref{lem:comm}(b), the fences of Table~\ref{tab:fence} and Theorem~\ref{thm:criterion} show that $\Hc_{e_k}$ is not coverable.
\end{proof}

\begin{table}[t]
\centering
\caption{The twelve intervals of $\mathcal{E}_\Hc$, with the value $k$ of $\NN$ just left of the right endpoint $e_k=\sqrt{8/(3k)}$. The right endpoint is non-coverable for $k\notin\{54,26,18\}$; $e_{18}=2/\sqrt{27}$ and $e_{54}=2/9$ are coverable, and $e_{26}=2/\sqrt{39}$ is open.}
\label{tab:Ehc}
\renewcommand{\arraystretch}{1.3}
\begin{tabular}{@{}rllr@{}}
\toprule
& interval of $\mathcal{E}_\Hc$ & decimal & $k$\\
\midrule
1 & $(0,\,\tfrac{2}{3\sqrt{14}})$ & $(0.000000,\,0.178174)$ & 84 \\
2 & $(\tfrac{1}{\sqrt{31}},\,\tfrac{2}{\sqrt{87}})$ & $(0.179605,\,0.214423)$ & 58 \\
3 & $(\tfrac{1}{\sqrt{21}},\,\tfrac{2}{9})$ & $(0.218218,\,0.222222)$ & 54 \\
4 & $(\tfrac{1}{\sqrt{19}},\,\tfrac{4}{\sqrt{222}})$ & $(0.229416,\,0.268462)$ & 37 \\
5 & $(\tfrac{8}{\sqrt{817}},\,\tfrac{2}{\sqrt{51}})$ & $(0.279885,\,0.280056)$ & 34 \\
6 & $(\tfrac{2}{7},\,\tfrac{2}{3\sqrt{5}})$ & $(0.285714,\,0.298142)$ & 30 \\
7 & $(\tfrac{29}{\sqrt{9331}},\,\tfrac{2}{\sqrt{42}})$ & $(0.300216,\,0.308607)$ & 28 \\
8 & $(\tfrac{3}{\sqrt{93}},\,\tfrac{4}{9\sqrt{2}})$ & $(0.311086,\,0.314270)$ & 27 \\
9 & $(\tfrac{12}{\sqrt{1443}},\,\tfrac{2}{\sqrt{39}})$ & $(0.315899,\,0.320256)$ & 26 \\
10 & $(\tfrac{2}{\sqrt{39}},\,\tfrac{4}{5\sqrt{6}})$ & $(0.320256,\,0.326599)$ & 25 \\
11 & $(\tfrac{1}{\sqrt{7}},\,\tfrac{2}{3\sqrt{3}})$ & $(0.377964,\,0.384900)$ & 18 \\
12 & $(\tfrac{8}{19},\,\tfrac{4}{3\sqrt{10}})$ & $(0.421053,\,0.421637)$ & 15 \\
\bottomrule
\end{tabular}
\end{table}

\begin{proof}[Proof of Corollary~\ref{cor:dmin}, honeycomb]
As for the square lattice, with $e_{84}$ in place of $d_{52}$ and $2/9\in\Cov_\Hc$ from Theorem~\ref{thm:new}.
\end{proof}

\begin{remark}
For $k\in\{54,26,18\}$ we have $2k\in\{108,52,36\}$, which are norms, so Lemma~\ref{lem:comm}(b) does not apply. For $k=18$ and $k=54$ commensurable hexagonal coverings indeed exist: the construction of \cite{PNF26} at $2/\sqrt{27}$ and the construction of Theorem~\ref{thm:new} at $2/9$. Fences exist for all three values, so by Lemma~\ref{lem:perfect} and the proof of Lemma~\ref{lem:line}, if $\Hc_{e_{26}}$ is coverable, then it has a covering whose centres form a hexagonal packing commensurable with $\Hc_{e_{26}}$, in which every disk contains exactly $26$ points.
\end{remark}

Since $\Hc_d\subset\mathcal{T}_d$, where $\mathcal{T}_d$ is the triangular lattice, non-coverability of $\Hc_d$ implies non-coverability of $\mathcal{T}_d$. All intervals of $\mathcal{E}_\Hc$ lie in the set of non-coverable spacings of $\mathcal{T}_d$ determined in \cite{KPS26}, so this gives no new information about the triangular lattice.

\begin{table}[t]
\centering
\caption{Fences for Lemma~\ref{lem:fence} at the borderline spacings. In every case $s=t\,(dv_1+dv_2)$ with the listed $t$. The last column is the maximum number of points in a closed unit disk centred anywhere on the fence, computed exactly (Remark~\ref{rem:float}).}
\label{tab:fence}
\renewcommand{\arraystretch}{1.15}
\begin{tabular}{@{}rccr@{\qquad}rccr@{}}
\toprule
\multicolumn{4}{c}{square lattice, $d=\sqrt{2\sqrt3/k}$} & \multicolumn{4}{c}{honeycomb, $d=\sqrt{8/(3k)}$}\\
\cmidrule(r){1-4}\cmidrule(l){5-8}
$k$ & $t$ & & max & $k$ & $t$ & $2k$ norm? & max\\
\midrule
52 & $1/80$ && 48 & 84 & $3/80$ & no & 82\\
44 & $1/80$ && 41 & 58 & $1/80$ & no & 54\\
37 & $9/80$ && 36 & 54 & $9/80$ & yes & 52\\
32 & $1/80$ && 30 & 37 & $1/80$ & no & 36\\
28 & $1/80$ && 26 & 34 & $5/80$ & no & 33\\
27 & $1/80$ && 26 & 30 & $1/80$ & no & 28\\
26 & $11/80$ && 25 & 28 & $3/80$ & no & 27\\
24 & $9/80$ && 23 & 27 & $1/80$ & no & 26\\
22 & $3/80$ && 21 & 26 & $1/80$ & yes & 24\\
21 & $21/80$ && 20 & 25 & $1/80$ & no & 24\\
17 & $19/80$ && 16 & 18 & $5/80$ & yes & 17\\
16 & $9/80$ && 15 & 15 & $1/80$ & no & 14\\
14 & $1/80$ && 13 & & & & \\
\bottomrule
\end{tabular}
\end{table}

\begin{figure}[t]
\centering
\includegraphics[width=0.95\linewidth]{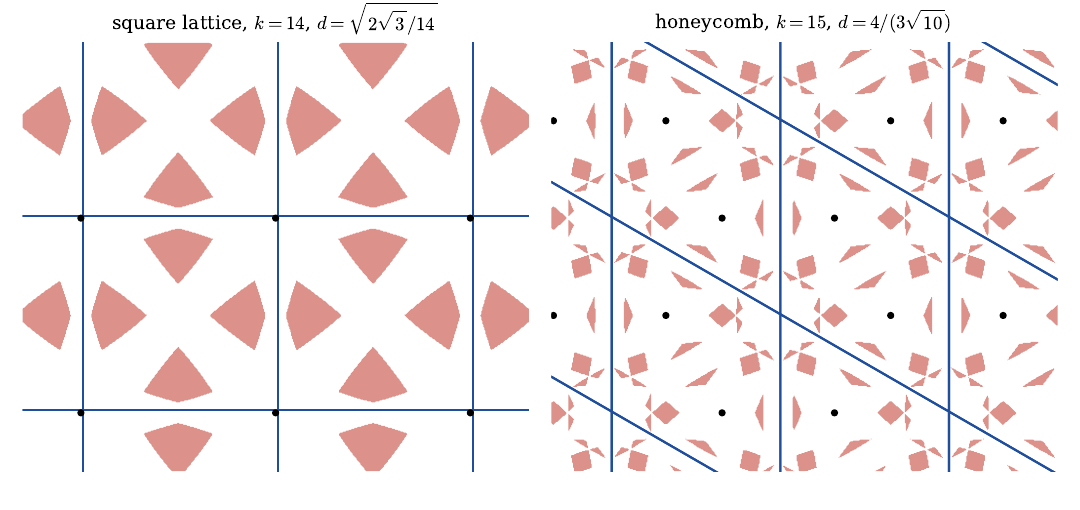}
\caption{The set $\Omega_k$ (red) of centres whose unit disk contains at least $k$ points (black), and the fence of Lemma~\ref{lem:fence} (blue), for the square lattice with $k=14$ and the honeycomb point set with $k=15$.}
\label{fig:fence}
\end{figure}

\section{New coverings}\label{sec:new}

A \emph{single-family lattice pattern} has centre set $d(c_0+M)$, where $M=\Z u+\Z v$ is a lattice spanned by two points $u,v$ of $\Z[i]$ (for the square lattice) or of $\Z[\omega]$ (for the honeycomb point set) and $c_0$ is a point with rational coordinates. Let $n_{\min}=\min_{0\ne m\in M}\Nm(m)$ and let $R$ be the covering radius of $\X$ by $c_0+M$, that is, $R=\sup_{x\in\X}\operatorname{dist}(x,c_0+M)$. The disks are non-overlapping if and only if $d\sqrt{n_{\min}}\ge2$, and they cover $\X_d$ if and only if $dR\le1$. Hence the pattern covers $\X_d$ for every
\[
d\in\bigl[2/\sqrt{n_{\min}},\,1/R\bigr].
\]
Since $\X$ is invariant under $3M\subseteq\Lambda_\X$ (for the honeycomb) or under $M\subseteq\Lambda_\X$ (for the square lattice), $R^2$ is a maximum over finitely many residue classes and can be computed exactly in rational arithmetic.

\begin{figure}[tbp]
\centering
\includegraphics[width=0.95\linewidth]{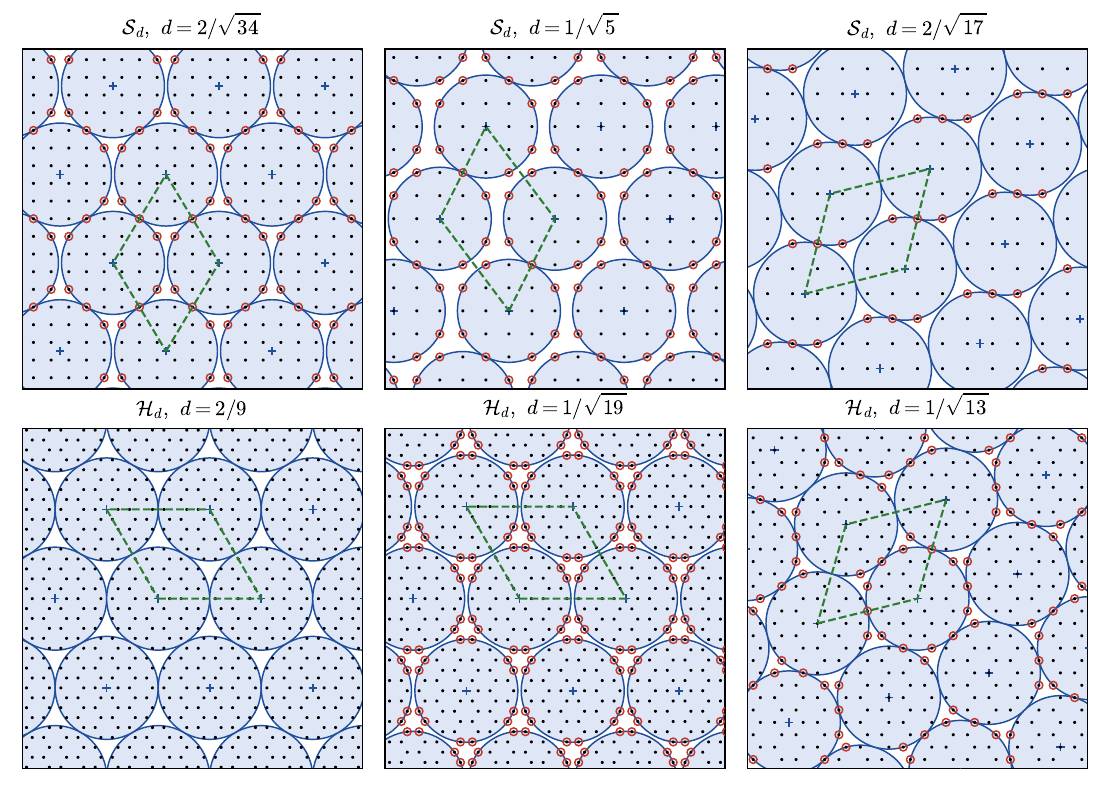}
\caption{The coverings of Table~\ref{tab:new}, in the window $[-3.3,3.3]^2$. Disks are non-overlapping closed unit disks centred at the points $+$ of $d(c_0+M)$; points of $\Sq_d$ or $\Hc_d$ at distance exactly $1$ from the nearest centre are circled in red, and all other points lie in the interior of a disk. The dashed parallelogram is a fundamental cell of $dM$. Top: the three square-lattice patterns, in each of which adjacent disks are tangent. Bottom: the honeycomb pattern with $u=9$, $v=9\omega$, $c_0=-3$ at both ends of its interval of spacings, where at $d=2/9$ the disks are tangent and no point lies on a boundary, and at $d=1/\sqrt{19}$ the disks are separated and every disk has points on its boundary; bottom right: the pattern $u=8+2\omega$, $v=6+8\omega$, $c_0=0$, restricted to $\Hc_d$.}
\label{fig:newcov}
\end{figure}

\begin{proof}[Proof of Theorem~\ref{thm:new}]
Table~\ref{tab:new} lists the patterns; Part~E of the program computes $n_{\min}$ and $R^2$ exactly. For the honeycomb pattern this can be checked by hand. Here $M=9\Z[\omega]$, so $n_{\min}=81$ and the disks are non-overlapping for $d\ge2/9$. There are $54$ points $x\in\Hc$ with $\Nm(x+3)\le19$. The difference of two of them has norm at most $4\cdot19=76<81$, so no two of them lie in the same class of $\Z[\omega]/9\Z[\omega]$. Since $9\Z[\omega]\subseteq\Lambda_\Hc$, membership in $\Hc$ depends only on the class modulo $9\Z[\omega]$, and $\Hc$ meets exactly $54$ of the $81$ classes (two thirds of them). Hence the $54$ points meet every class that contains points of $\Hc$, every point of $\Hc$ lies within distance $\sqrt{19}$ of $-3+M$, and the disks cover $\Hc_d$ for $d\le1/\sqrt{19}$. The last pattern covers every point of $\Z[\omega]$, so it covers the triangular lattice $\mathcal T_d$ at $d=1/\sqrt{13}$ (this covering is also given in \cite[Theorem~C]{KPS26}) and, a fortiori, its subset $\Hc_d$. The coverings are shown in Figure~\ref{fig:newcov}. The remaining square-lattice spacings of Theorem~\ref{thm:new} are obtained after Lemma~\ref{lem:sub} below.
\end{proof}

\begin{table}[t]
\centering
\caption{Single-family lattice patterns proving Theorem~\ref{thm:new}. Coordinates are Gaussian ($a+bi$) for the square lattice and Eisenstein ($a+b\omega$) for the honeycomb point set. The last pattern covers all of $\Z[\omega]$. The covering at $2/\sqrt{17}$ also follows from the one at $2/\sqrt{34}$ by Lemma~\ref{lem:sub}(a).}
\label{tab:new}
\renewcommand{\arraystretch}{1.25}
\begin{tabular}{@{}lllllll@{}}
\toprule
point set & $u$ & $v$ & $c_0$ & $n_{\min}$ & $R^2$ & coverable for\\
\midrule
$\Sq$ & $3+5i$ & $-3+5i$ & $-\frac32+\frac52i$ & $34$ & $17/2$ & $d=2/\sqrt{34}$\\
$\Sq$ & $2+4i$ & $-3+4i$ & $0$ & $20$ & $5$ & $d=1/\sqrt5$\\
$\Sq$ & $4+i$ & $1+4i$ & $\frac12+2i$ & $17$ & $17/4$ & $d=2/\sqrt{17}$\\
$\Hc$ & $9$ & $9\omega$ & $-3$ & $81$ & $19$ & $2/9\le d\le1/\sqrt{19}$\\
$\Z[\omega]\supset\Hc$ & $8+2\omega$ & $6+8\omega$ & $0$ & $52$ & $13$ & $d=1/\sqrt{13}$\\
\bottomrule
\end{tabular}
\end{table}

\begin{remark}
In the honeycomb pattern of Table~\ref{tab:new}, each disk is centred at a point of the triangular lattice that is not in $\Hc$ (a centre of a hexagon of the honeycomb) and covers the $54$ points of $\Hc$ within distance $\sqrt{19}\,d$. At $d=2/9$ the centres form a hexagonal packing commensurable with $\Hc_d$, as Lemma~\ref{lem:perfect} requires at the borderline spacing $e_{54}$. The search that found these patterns was limited to centre offsets $c_0$ in a grid of mesh $1/6$ (refined numerically) and is not claimed to be exhaustive.
\end{remark}

\subsection{Sublattices}

Each of the three point sets contains scaled and rotated copies of itself, and the honeycomb point set and the triangular lattice contain copies of each other. This relates different spacings.

\begin{lemma}\label{lem:sub}
\begin{enumerate}[label=(\alph*),nosep]
\item Let $\alpha\in\Z[i]\setminus\{0\}$ and $q=\Nm(\alpha)$. Then $\alpha\Sq_d$, a rotated copy of $\Sq_{d\sqrt q}$, is contained in $\Sq_d$. Hence, if $\Sq_d$ is coverable, so is $\Sq_{d\sqrt q}$; if $\Sq_{d\sqrt q}$ is not coverable, neither is $\Sq_d$.
\item Let $\alpha\in\Z[\omega]\setminus\{0\}$ and $q=\Nm(\alpha)$. Then a translate of $\alpha\Hc_d$, which is a rotated copy of $\Hc_{d\sqrt q}$, is contained in $\Hc_d$, and the same two conclusions hold for $\Hc$.
\item $\Hc_d\subseteq\mathcal T_d$, and a translate of $(1-\omega)\mathcal T_d$, a rotated copy of $\mathcal T_{d\sqrt3}$, is contained in $\Hc_d$. Hence coverability of $\mathcal T_d$ implies that of $\Hc_d$, and coverability of $\Hc_d$ implies that of $\mathcal T_{d\sqrt3}$.
\end{enumerate}
\end{lemma}

\begin{proof}
(a) holds because $\alpha\Z[i]\subseteq\Z[i]$ and multiplication by $\alpha$ is a rotation composed with a scaling by $\sqrt q$; coverability is invariant under rotations and translations and inherited by subsets. (c) $\Hc\subseteq\Z[\omega]$ by definition, and $1+(1-\omega)\Z[\omega]=1+\Lambda_\Hc\subseteq\Hc$, because $a+b\omega\in1+\Lambda_\Hc$ implies $a+b\equiv1\pmod3$. (b) Recall $\Hc=(1+\Lambda_\Hc)\cup(-1+\Lambda_\Hc)$ and $\alpha\Lambda_\Hc\subseteq\Lambda_\Hc$. If $3\nmid q$, then $1-\omega\nmid\alpha$, so $\alpha\equiv\varepsilon\pmod{\Lambda_\Hc}$ for some $\varepsilon\in\{1,-1\}$, because $\Z[\omega]/\Lambda_\Hc$ has the three classes of $0,1,-1$. Then $\alpha(\pm1+\Lambda_\Hc)\subseteq\pm\varepsilon+\Lambda_\Hc\subseteq\Hc$, so $\alpha\Hc\subseteq\Hc$. If $3\mid q$, then $(1-\omega)\mid\alpha$, so $\alpha\Hc\subseteq\Lambda_\Hc$ and $1+\alpha\Hc\subseteq\Hc$ by (c).
\end{proof}

By Lemma~\ref{lem:sub}(a) with $\alpha=1+i$, the coverable interval $[1/2,6/\sqrt{130}]$ of Theorem~\ref{thm:PNF} gives $[1/\sqrt2,6/\sqrt{65}]$, and with $\alpha=1+2i$ the coverable spacing $2/\sqrt{34}$ gives $\sqrt5\cdot2/\sqrt{34}=\sqrt{10/17}$. This completes the proof of Theorem~\ref{thm:new}.

Applying Lemma~\ref{lem:sub} to Theorem~\ref{thm:PNF}, Theorems~\ref{thm:sq}--\ref{thm:new} and to the results of \cite{KPS26} for the triangular lattice gives exactly the three square-lattice statements $2/\sqrt{17}=\sqrt2\cdot2/\sqrt{34}$, $\sqrt{10/17}=\sqrt5\cdot2/\sqrt{34}$ and $[1/\sqrt2,6/\sqrt{65}]=\sqrt2\cdot[1/2,6/\sqrt{130}]$ of Theorem~\ref{thm:new}, and nothing new for the honeycomb point set; Part~F of the program checks this in exact arithmetic. Only small $q$ matter. For the square lattice, every known non-coverable spacing is at most $d_{14}$ and every open interval lies above $d_{52}$, so $d/\sqrt q$ can meet an open interval only for $q=2$; every known coverable spacing is at least $2/\sqrt{34}$, so $d\sqrt q$ can meet an open interval only for $q\in\{2,4,5\}$. For the honeycomb point set the same comparison leaves $q\in\{3,4\}$ and $q\in\{3,4,7\}$, respectively; Part~F also checks the transfers between $\Hc_d$ and the triangular lattice in Lemma~\ref{lem:sub}(c). The lemma also links open intervals. For example, $\Hc_d$ contains a translated and rotated copy of $\mathcal T_{3d}$ by Lemma~\ref{lem:sub}(b),(c), and $3\cdot(e_{84},1/\sqrt{31}]=(2/\sqrt{14},3/\sqrt{31}]$; hence, if $\mathcal T_d$ is non-coverable for every $d\in(2/\sqrt{14},3/\sqrt{31}]$, then $\Hc_d$ is non-coverable for every $d\in(e_{84},1/\sqrt{31}]$ and $d_{\min}^\Hc>2/\sqrt{87}$. The same hypothesis on the triangular lattice controls the window containing its smallest coverable spacing \cite[Proposition~9.2]{KPS26}.

\subsection{Tight and loose coverings}

Proposition~\ref{prop:slack}(c) applies wherever a coverable spacing is followed or preceded by non-coverable spacings. For the square lattice, the spacings $2/\sqrt{34}$, $1/\sqrt5$ and $2/\sqrt{17}$ of Theorem~\ref{thm:new} are left endpoints of intervals of $\mathcal E_\Sq$, so every covering at these spacings has $m(C)=1$; the coverings of Figure~\ref{fig:newcov} indeed have points on the boundaries of their disks. For the honeycomb point set, every covering has $m(C)=1$ at $1/\sqrt7$ and at $1/\sqrt{19}$, which are followed by the intervals $(1/\sqrt7,e_{18})$ and $(1/\sqrt{19},e_{37})$ of $\mathcal E_\Hc$, and every covering has $s(C)=2$ at $e_{18}=2/\sqrt{27}$ and $e_{54}=2/9$, which are preceded by intervals of $\mathcal E_\Hc$. At $2/9$ this agrees with Lemma~\ref{lem:perfect}.

In the other direction, every spacing in the interior of the known coverable sets has a covering with $m(C)<1$, so in particular one with no point on the boundary of a disk. For the constructions of Table~\ref{tab:new} and all constructions of \cite{PNF26} for $\Hc_d$ \cite[Section~4, Table~2]{PNF26}, and for those for $\Sq_d$ except the zipper pattern \cite[Section~5, Table~3]{PNF26}, the disks are in fixed relative position, so all distances scale linearly with $d$, and each construction is valid on an interval $[d_1,d_2]$ with $s=2$ at $d_1$ and $m=1$ at $d_2$; at every junction of two such intervals in Theorem~\ref{thm:PNF}, the spacing is the left endpoint of one of them. The zipper pattern of \cite[Section~5.3, Figure~14]{PNF26}, valid for $1/2\le d\le6/\sqrt{130}$, uses centres that depend on $d$. On $[1/2,4/\sqrt{61}]$ the centres are in fixed relative position, as above. On $[4/\sqrt{61},6/\sqrt{130}]$ two lattice points of each motif lie on the boundary of its disk for every $d$. Its non-overlap constraints are inactive there, so $s(C)>2$, and Proposition~\ref{prop:slack}(b) applied at a slightly larger spacing gives a covering with $m<1$ and $s>2$ at every $d\in[4/\sqrt{61},6/\sqrt{130})$. Finally, coverings obtained from Lemma~\ref{lem:sub} are restrictions of coverings at smaller spacings and inherit their values of $m$ and $s$. Whether every interior point of $\Cov_\X$ has a covering with no point on the boundary of a disk is open (Section~\ref{sec:open}).

\section{Open problems}\label{sec:open}

\begin{enumerate}[label=(\arabic*)]
\item Is $\Hc_d$ coverable at $d=2/\sqrt{39}$? By Lemma~\ref{lem:perfect}, if it is, then it has a covering whose centres form a hexagonal packing, commensurable with $\Hc_d$, in which every disk contains exactly $26$ points.
\item Are the gaps $(6/\sqrt{130},2/\sqrt{13})$, $(2/\sqrt{10},4/\sqrt{37})$ and $(6/\sqrt{65},4/\sqrt{26})\setminus\{\sqrt{10/17}\}$ for the square lattice, and $(5/\sqrt{91},2/\sqrt{13})$ and $(3/\sqrt{21},4/\sqrt{37})$ for the honeycomb point set, coverable? On these gaps the density ratio stays above $1.02$ (its infimum is about $1.025$, $1.039$, $1.119$ for the square lattice and $1.087$, $1.125$ for the honeycomb point set), so the density argument does not apply.
\item What are $d_{\min}^\Sq$ and $d_{\min}^\Hc$?
\item Let $d$ be an interior point of $\Cov_\X$. Does $\X_d$ have a covering in which no point lies on the boundary of a disk? By Proposition~\ref{prop:slack}, the answer is yes at $d$ as soon as, for some $d'>d$, the set $\X_{d'}$ has a covering $C$ with $s(C)>2d'/d$, or, for some $d'<d$, it has a covering $C$ with $m(C)<d'/d$.
\end{enumerate}

\section*{Acknowledgements}
The research of this paper was carried out while the authors were participating in the Geometry and Combinatorics Boot Camp 2026, in Lamphun, Thailand. We are grateful for the hospitality of the Department of Mathematics, Faculty of Science, Chiang Mai University.

\section*{Declaration of generative AI use}
During the preparation of this work the authors used Claude (Anthropic) to assist with extending the ideas of the companion paper on the triangular lattice \cite{KPS26} to the square lattice and the honeycomb point set, computer searches, writing and checking the verification program, and drafting and editing the text. The authors checked all proofs, computations and text, and take full responsibility for the content of this paper.


\appendix
\section{The step functions \texorpdfstring{$\NN$}{N}}\label{app:tables}

The last column gives the part of each block that belongs to $\mathcal E$; the column $\NN b^2/(2\sqrt3\rho_1)$ is the density ratio at the right end of the block.

{\small
\renewcommand{\arraystretch}{1.2}
\begin{longtable}{@{}rllllrl@{}}
\caption{Blocks of constancy of $\NN$ for the square lattice on $[1/5,2]$.}\label{tab:blocksq}\\
\toprule
\# & block $(a,b]$ & $a$ & $b$ & $\NN$ & ratio & part in $\mathcal{E}_\Sq$\\
\midrule
\endfirsthead
\toprule
\# & block $(a,b]$ & $a$ & $b$ & $\NN$ & ratio & part in $\mathcal{E}_\Sq$\\
\midrule
\endhead
\bottomrule
\endlastfoot
1 & $[\tfrac{1}{5},\ \tfrac{14}{\sqrt{4810}}]$ & $0.200000$ & $0.201862$ & 81 & $0.9528$ & $(\tfrac{1}{5},\,\tfrac{14}{\sqrt{4810}}]$ \\
2 & $(\tfrac{14}{\sqrt{4810}},\ \tfrac{62}{85\sqrt{13}}]$ & $0.201862$ & $0.202302$ & 80 & $0.9452$ & $(\tfrac{14}{\sqrt{4810}},\,\tfrac{62}{85\sqrt{13}}]$ \\
3 & $(\tfrac{62}{85\sqrt{13}},\ \tfrac{86}{5\sqrt{7178}}]$ & $0.202302$ & $0.203014$ & 79 & $0.9399$ & $(\tfrac{62}{85\sqrt{13}},\,\tfrac{86}{5\sqrt{7178}}]$ \\
4 & $(\tfrac{86}{5\sqrt{7178}},\ \tfrac{6}{5\sqrt{34}}]$ & $0.203014$ & $0.205798$ & 78 & $0.9536$ & $(\tfrac{86}{5\sqrt{7178}},\,\tfrac{6}{5\sqrt{34}}]$ \\
5 & $(\tfrac{6}{5\sqrt{34}},\ \tfrac{26}{89\sqrt{2}}]$ & $0.205798$ & $0.206571$ & 77 & $0.9485$ & $(\tfrac{6}{5\sqrt{34}},\,\tfrac{26}{89\sqrt{2}}]$ \\
6 & $(\tfrac{26}{89\sqrt{2}},\ \tfrac{2}{3\sqrt{10}}]$ & $0.206571$ & $0.210819$ & 76 & $0.9751$ & $(\tfrac{26}{89\sqrt{2}},\,\tfrac{2}{3\sqrt{10}}]$ \\
7 & $(\tfrac{2}{3\sqrt{10}},\ \tfrac{2}{\sqrt{89}}]$ & $0.210819$ & $0.212000$ & 74 & $0.9601$ & $(\tfrac{2}{3\sqrt{10}},\,\tfrac{2}{\sqrt{89}}]$ \\
8 & $(\tfrac{2}{\sqrt{89}},\ \tfrac{3}{10\sqrt{2}}]$ & $0.212000$ & $0.212132$ & 73 & $0.9483$ & $(\tfrac{2}{\sqrt{89}},\,\tfrac{3}{10\sqrt{2}}]$ \\
9 & $(\tfrac{3}{10\sqrt{2}},\ \tfrac{14}{5\sqrt{170}}]$ & $0.212132$ & $0.214750$ & 72 & $0.9585$ & $(\tfrac{3}{10\sqrt{2}},\,\tfrac{14}{5\sqrt{170}}]$ \\
10 & $(\tfrac{14}{5\sqrt{170}},\ \tfrac{34}{\sqrt{25010}}]$ & $0.214750$ & $0.214992$ & 71 & $0.9474$ & $(\tfrac{14}{5\sqrt{170}},\,\tfrac{34}{\sqrt{25010}}]$ \\
11 & $(\tfrac{34}{\sqrt{25010}},\ \tfrac{2}{\sqrt{85}}]$ & $0.214992$ & $0.216930$ & 70 & $0.9509$ & $(\tfrac{34}{\sqrt{25010}},\,\tfrac{2}{\sqrt{85}}]$ \\
12 & $(\tfrac{2}{\sqrt{85}},\ \tfrac{1}{2\sqrt{5}}]$ & $0.216930$ & $0.223607$ & 69 & $0.9959$ & $(\tfrac{2}{\sqrt{85}},\,\tfrac{1}{2\sqrt{5}}]$ \\
13 & $(\tfrac{1}{2\sqrt{5}},\ \tfrac{5}{\sqrt{481}}]$ & $0.223607$ & $0.227980$ & 65 & $0.9753$ & $(\tfrac{1}{2\sqrt{5}},\,\tfrac{5}{\sqrt{481}}]$ \\
14 & $(\tfrac{5}{\sqrt{481}},\ \tfrac{22}{5\sqrt{370}}]$ & $0.227980$ & $0.228745$ & 64 & $0.9667$ & $(\tfrac{5}{\sqrt{481}},\,\tfrac{22}{5\sqrt{370}}]$ \\
15 & $(\tfrac{22}{5\sqrt{370}},\ \tfrac{66}{\sqrt{81770}}]$ & $0.228745$ & $0.230806$ & 63 & $0.9688$ & $(\tfrac{22}{5\sqrt{370}},\,\tfrac{66}{\sqrt{81770}}]$ \\
16 & $(\tfrac{66}{\sqrt{81770}},\ \tfrac{22}{\sqrt{8906}}]$ & $0.230806$ & $0.233121$ & 62 & $0.9727$ & $(\tfrac{66}{\sqrt{81770}},\,\tfrac{22}{\sqrt{8906}}]$ \\
17 & $(\tfrac{22}{\sqrt{8906}},\ \tfrac{1}{3\sqrt{2}}]$ & $0.233121$ & $0.235702$ & 61 & $0.9783$ & $(\tfrac{22}{\sqrt{8906}},\,\tfrac{1}{3\sqrt{2}}]$ \\
18 & $(\tfrac{1}{3\sqrt{2}},\ \tfrac{14}{\sqrt{3445}}]$ & $0.235702$ & $0.238525$ & 59 & $0.9690$ & $(\tfrac{1}{3\sqrt{2}},\,\tfrac{14}{\sqrt{3445}}]$ \\
19 & $(\tfrac{14}{\sqrt{3445}},\ \tfrac{8}{\sqrt{1105}}]$ & $0.238525$ & $0.240663$ & 58 & $0.9697$ & $(\tfrac{14}{\sqrt{3445}},\,\tfrac{8}{\sqrt{1105}}]$ \\
20 & $(\tfrac{8}{\sqrt{1105}},\ \tfrac{1}{\sqrt{17}}]$ & $0.240663$ & $0.242536$ & 57 & $0.9679$ & $(\tfrac{8}{\sqrt{1105}},\,\tfrac{1}{\sqrt{17}}]$ \\
21 & $(\tfrac{1}{\sqrt{17}},\ \tfrac{2}{\sqrt{65}}]$ & $0.242536$ & $0.248069$ & 56 & $0.9948$ & $(\tfrac{1}{\sqrt{17}},\,\tfrac{2}{\sqrt{65}}]$ \\
22 & $(\tfrac{2}{\sqrt{65}},\ \tfrac{2}{\sqrt{58}}]$ & $0.248069$ & $0.262613$ & 52 & $1.0352$ & $(\tfrac{2}{\sqrt{65}},\,\sqrt{2\sqrt3/52})$ \\
23 & $(\tfrac{2}{\sqrt{58}},\ \tfrac{10}{\sqrt{1394}}]$ & $0.262613$ & $0.267836$ & 48 & $0.9940$ & $(\tfrac{2}{\sqrt{58}},\,\tfrac{10}{\sqrt{1394}}]$ \\
24 & $(\tfrac{10}{\sqrt{1394}},\ \tfrac{66}{5\sqrt{2378}}]$ & $0.267836$ & $0.270687$ & 47 & $0.9941$ & $(\tfrac{10}{\sqrt{1394}},\,\tfrac{66}{5\sqrt{2378}}]$ \\
25 & $(\tfrac{66}{5\sqrt{2378}},\ \tfrac{5}{13\sqrt{2}}]$ & $0.270687$ & $0.271964$ & 46 & $0.9822$ & $(\tfrac{66}{5\sqrt{2378}},\,\tfrac{5}{13\sqrt{2}}]$ \\
26 & $(\tfrac{5}{13\sqrt{2}},\ \tfrac{1}{\sqrt{13}}]$ & $0.271964$ & $0.277350$ & 45 & $0.9993$ & $(\tfrac{5}{13\sqrt{2}},\,\tfrac{1}{\sqrt{13}}]$ \\
27 & $(\tfrac{1}{\sqrt{13}},\ \tfrac{2}{5\sqrt{2}}]$ & $0.277350$ & $0.282843$ & 44 & $1.0161$ & $(\tfrac{1}{\sqrt{13}},\,\sqrt{2\sqrt3/44})$ \\
28 & $(\tfrac{2}{5\sqrt{2}},\ \tfrac{10}{\sqrt{1189}}]$ & $0.282843$ & $0.290007$ & 41 & $0.9954$ & $(\tfrac{2}{5\sqrt{2}},\,\tfrac{10}{\sqrt{1189}}]$ \\
29 & $(\tfrac{10}{\sqrt{1189}},\ \tfrac{58}{\sqrt{39442}}]$ & $0.290007$ & $0.292044$ & 40 & $0.9848$ & $(\tfrac{10}{\sqrt{1189}},\,\tfrac{58}{\sqrt{39442}}]$ \\
30 & $(\tfrac{58}{\sqrt{39442}},\ \tfrac{6}{\sqrt{410}}]$ & $0.292044$ & $0.296319$ & 39 & $0.9885$ & $(\tfrac{58}{\sqrt{39442}},\,\tfrac{6}{\sqrt{410}}]$ \\
31 & $(\tfrac{6}{\sqrt{410}},\ \tfrac{2}{3\sqrt{5}}]$ & $0.296319$ & $0.298142$ & 38 & $0.9751$ & $(\tfrac{6}{\sqrt{410}},\,\tfrac{2}{3\sqrt{5}}]$ \\
32 & $(\tfrac{2}{3\sqrt{5}},\ \tfrac{1}{\sqrt{10}}]$ & $0.298142$ & $0.316228$ & 37 & $1.0681$ & $(\tfrac{2}{3\sqrt{5}},\,\sqrt{2\sqrt3/37})$ \\
33 & $(\tfrac{1}{\sqrt{10}},\ \tfrac{10}{\sqrt{962}}]$ & $0.316228$ & $0.322413$ & 33 & $0.9903$ & $(\tfrac{1}{\sqrt{10}},\,\tfrac{10}{\sqrt{962}}]$ \\
34 & $(\tfrac{10}{\sqrt{962}},\ \tfrac{2}{\sqrt{34}}]$ & $0.322413$ & $0.342997$ & 32 & $1.0868$ & $(\tfrac{10}{\sqrt{962}},\,\sqrt{2\sqrt3/32})$ \\
35 & $(\tfrac{2}{\sqrt{34}},\ \tfrac{4}{5\sqrt{5}}]$ & $0.342997$ & $0.357771$ & 28 & $1.0346$ & $(\tfrac{2}{\sqrt{34}},\,\sqrt{2\sqrt3/28})$ \\
36 & $(\tfrac{4}{5\sqrt{5}},\ \tfrac{38}{5\sqrt{442}}]$ & $0.357771$ & $0.361495$ & 27 & $1.0185$ & $(\tfrac{4}{5\sqrt{5}},\,\sqrt{2\sqrt3/27})$ \\
37 & $(\tfrac{38}{5\sqrt{442}},\ \tfrac{2}{\sqrt{29}}]$ & $0.361495$ & $0.371391$ & 26 & $1.0352$ & $(\tfrac{38}{5\sqrt{442}},\,\sqrt{2\sqrt3/26})$ \\
38 & $(\tfrac{2}{\sqrt{29}},\ \tfrac{2}{\sqrt{26}}]$ & $0.371391$ & $0.392232$ & 24 & $1.0659$ & $(\tfrac{2}{\sqrt{29}},\,\sqrt{2\sqrt3/24})$ \\
39 & $(\tfrac{2}{\sqrt{26}},\ \tfrac{2}{5}]$ & $0.392232$ & $0.400000$ & 22 & $1.0161$ & $(\tfrac{2}{\sqrt{26}},\,\sqrt{2\sqrt3/22})$ \\
40 & $(\tfrac{2}{5},\ \tfrac{1}{\sqrt{5}}]$ & $0.400000$ & $0.447214$ & 21 & $1.2124$ & $(\tfrac{2}{5},\,\sqrt{2\sqrt3/21})$ \\
41 & $(\tfrac{1}{\sqrt{5}},\ \tfrac{6}{13}]$ & $0.447214$ & $0.461538$ & 17 & $1.0454$ & $(\tfrac{1}{\sqrt{5}},\,\sqrt{2\sqrt3/17})$ \\
42 & $(\tfrac{6}{13},\ \tfrac{2}{\sqrt{17}}]$ & $0.461538$ & $0.485071$ & 16 & $1.0868$ & $(\tfrac{6}{13},\,\sqrt{2\sqrt3/16})$ \\
43 & $(\tfrac{2}{\sqrt{17}},\ \tfrac{6}{\sqrt{130}}]$ & $0.485071$ & $0.526235$ & 14 & $1.1192$ & $(\tfrac{2}{\sqrt{17}},\,\sqrt{2\sqrt3/14})$ \\
44 & $(\tfrac{6}{\sqrt{130}},\ \tfrac{10}{13\sqrt{2}}]$ & $0.526235$ & $0.543928$ & 13 & $1.1103$ & --- \\
45 & $(\tfrac{10}{13\sqrt{2}},\ \tfrac{2}{\sqrt{10}}]$ & $0.543928$ & $0.632456$ & 12 & $1.3856$ & --- \\
46 & $(\tfrac{2}{\sqrt{10}},\ \tfrac{1}{\sqrt{2}}]$ & $0.632456$ & $0.707107$ & 9 & $1.2990$ & --- \\
47 & $(\tfrac{1}{\sqrt{2}},\ \tfrac{4}{5}]$ & $0.707107$ & $0.800000$ & 7 & $1.2933$ & --- \\
48 & $(\tfrac{4}{5},\ \tfrac{2}{\sqrt{5}}]$ & $0.800000$ & $0.894427$ & 6 & $1.3856$ & --- \\
49 & $(\tfrac{2}{\sqrt{5}},\ 1]$ & $0.894427$ & $1.000000$ & 5 & $1.4434$ & --- \\
50 & $(1,\ \tfrac{2}{\sqrt{2}}]$ & $1.000000$ & $1.414214$ & 4 & $2.3094$ & --- \\
51 & $(\tfrac{2}{\sqrt{2}},\ 2]$ & $1.414214$ & $2.000000$ & 2 & $2.3094$ & --- \\
\end{longtable}

\begin{longtable}{@{}rllllrl@{}}
\caption{Blocks of constancy of $\NN$ for the honeycomb point set on $[1/7,2]$.}\label{tab:blockhc}\\
\toprule
\# & block $(a,b]$ & $a$ & $b$ & $\NN$ & ratio & part in $\mathcal{E}_\Hc$\\
\midrule
\endfirsthead
\toprule
\# & block $(a,b]$ & $a$ & $b$ & $\NN$ & ratio & part in $\mathcal{E}_\Hc$\\
\midrule
\endhead
\bottomrule
\endlastfoot
1 & $[\tfrac{1}{7},\ \tfrac{73}{\sqrt{261079}}]$ & $0.142857$ & $0.142869$ & 122 & $0.9338$ & $(\tfrac{1}{7},\,\tfrac{73}{\sqrt{261079}}]$ \\
2 & $(\tfrac{73}{\sqrt{261079}},\ \tfrac{5}{\sqrt{1209}}]$ & $0.142869$ & $0.143799$ & 121 & $0.9383$ & $(\tfrac{73}{\sqrt{261079}},\,\tfrac{5}{\sqrt{1209}}]$ \\
3 & $(\tfrac{5}{\sqrt{1209}},\ \tfrac{13}{\sqrt{8113}}]$ & $0.143799$ & $0.144329$ & 120 & $0.9374$ & $(\tfrac{5}{\sqrt{1209}},\,\tfrac{13}{\sqrt{8113}}]$ \\
4 & $(\tfrac{13}{\sqrt{8113}},\ \tfrac{139}{13\sqrt{5461}}]$ & $0.144329$ & $0.144689$ & 119 & $0.9342$ & $(\tfrac{13}{\sqrt{8113}},\,\tfrac{139}{13\sqrt{5461}}]$ \\
5 & $(\tfrac{139}{13\sqrt{5461}},\ \tfrac{5}{13\sqrt{7}}]$ & $0.144689$ & $0.145371$ & 118 & $0.9351$ & $(\tfrac{139}{13\sqrt{5461}},\,\tfrac{5}{13\sqrt{7}}]$ \\
6 & $(\tfrac{5}{13\sqrt{7}},\ \tfrac{2}{3\sqrt{21}}]$ & $0.145371$ & $0.145479$ & 117 & $0.9286$ & $(\tfrac{5}{13\sqrt{7}},\,\tfrac{2}{3\sqrt{21}}]$ \\
7 & $(\tfrac{2}{3\sqrt{21}},\ \tfrac{5}{\sqrt{1147}}]$ & $0.145479$ & $0.147635$ & 116 & $0.9481$ & $(\tfrac{2}{3\sqrt{21}},\,\tfrac{5}{\sqrt{1147}}]$ \\
8 & $(\tfrac{5}{\sqrt{1147}},\ \tfrac{68}{7\sqrt{4309}}]$ & $0.147635$ & $0.147987$ & 115 & $0.9444$ & $(\tfrac{5}{\sqrt{1147}},\,\tfrac{68}{7\sqrt{4309}}]$ \\
9 & $(\tfrac{68}{7\sqrt{4309}},\ \tfrac{2}{\sqrt{181}}]$ & $0.147987$ & $0.148659$ & 114 & $0.9448$ & $(\tfrac{68}{7\sqrt{4309}},\,\tfrac{2}{\sqrt{181}}]$ \\
10 & $(\tfrac{2}{\sqrt{181}},\ \tfrac{131}{\sqrt{772597}}]$ & $0.148659$ & $0.149037$ & 113 & $0.9412$ & $(\tfrac{2}{\sqrt{181}},\,\tfrac{131}{\sqrt{772597}}]$ \\
11 & $(\tfrac{131}{\sqrt{772597}},\ \tfrac{3}{\sqrt{399}}]$ & $0.149037$ & $0.150188$ & 112 & $0.9474$ & $(\tfrac{131}{\sqrt{772597}},\,\tfrac{3}{\sqrt{399}}]$ \\
12 & $(\tfrac{3}{\sqrt{399}},\ \tfrac{107}{91\sqrt{61}}]$ & $0.150188$ & $0.150549$ & 111 & $0.9434$ & $(\tfrac{3}{\sqrt{399}},\,\tfrac{107}{91\sqrt{61}}]$ \\
13 & $(\tfrac{107}{91\sqrt{61}},\ \tfrac{2}{5\sqrt{7}}]$ & $0.150549$ & $0.151186$ & 110 & $0.9429$ & $(\tfrac{107}{91\sqrt{61}},\,\tfrac{2}{5\sqrt{7}}]$ \\
14 & $(\tfrac{2}{5\sqrt{7}},\ \tfrac{1}{\sqrt{43}}]$ & $0.151186$ & $0.152499$ & 109 & $0.9506$ & $(\tfrac{2}{5\sqrt{7}},\,\tfrac{1}{\sqrt{43}}]$ \\
15 & $(\tfrac{1}{\sqrt{43}},\ \tfrac{9}{\sqrt{3477}}]$ & $0.152499$ & $0.152630$ & 108 & $0.9435$ & $(\tfrac{1}{\sqrt{43}},\,\tfrac{9}{\sqrt{3477}}]$ \\
16 & $(\tfrac{9}{\sqrt{3477}},\ \tfrac{17}{7\sqrt{247}}]$ & $0.152630$ & $0.154526$ & 107 & $0.9581$ & $(\tfrac{9}{\sqrt{3477}},\,\tfrac{17}{7\sqrt{247}}]$ \\
17 & $(\tfrac{17}{7\sqrt{247}},\ \tfrac{31}{7\sqrt{817}}]$ & $0.154526$ & $0.154936$ & 106 & $0.9542$ & $(\tfrac{17}{7\sqrt{247}},\,\tfrac{31}{7\sqrt{817}}]$ \\
18 & $(\tfrac{31}{7\sqrt{817}},\ \tfrac{25}{\sqrt{25591}}]$ & $0.154936$ & $0.156277$ & 105 & $0.9616$ & $(\tfrac{31}{7\sqrt{817}},\,\tfrac{25}{\sqrt{25591}}]$ \\
19 & $(\tfrac{25}{\sqrt{25591}},\ \tfrac{83}{\sqrt{281827}}]$ & $0.156277$ & $0.156346$ & 104 & $0.9533$ & $(\tfrac{25}{\sqrt{25591}},\,\tfrac{83}{\sqrt{281827}}]$ \\
20 & $(\tfrac{83}{\sqrt{281827}},\ \tfrac{1}{\sqrt{39}}]$ & $0.156346$ & $0.160128$ & 103 & $0.9904$ & $(\tfrac{83}{\sqrt{281827}},\,\tfrac{1}{\sqrt{39}}]$ \\
21 & $(\tfrac{1}{\sqrt{39}},\ \tfrac{13}{7\sqrt{133}}]$ & $0.160128$ & $0.161035$ & 97 & $0.9433$ & $(\tfrac{1}{\sqrt{39}},\,\tfrac{13}{7\sqrt{133}}]$ \\
22 & $(\tfrac{13}{7\sqrt{133}},\ \tfrac{1}{\sqrt{37}}]$ & $0.161035$ & $0.164399$ & 96 & $0.9730$ & $(\tfrac{13}{7\sqrt{133}},\,\tfrac{1}{\sqrt{37}}]$ \\
23 & $(\tfrac{1}{\sqrt{37}},\ \tfrac{2}{7\sqrt{3}}]$ & $0.164399$ & $0.164957$ & 92 & $0.9388$ & $(\tfrac{1}{\sqrt{37}},\,\tfrac{2}{7\sqrt{3}}]$ \\
24 & $(\tfrac{2}{7\sqrt{3}},\ \tfrac{9}{\sqrt{2847}}]$ & $0.164957$ & $0.168674$ & 90 & $0.9602$ & $(\tfrac{2}{7\sqrt{3}},\,\tfrac{9}{\sqrt{2847}}]$ \\
25 & $(\tfrac{9}{\sqrt{2847}},\ \tfrac{17}{\sqrt{10101}}]$ & $0.168674$ & $0.169148$ & 89 & $0.9549$ & $(\tfrac{9}{\sqrt{2847}},\,\tfrac{17}{\sqrt{10101}}]$ \\
26 & $(\tfrac{17}{\sqrt{10101}},\ \tfrac{7}{\sqrt{1677}}]$ & $0.169148$ & $0.170935$ & 88 & $0.9642$ & $(\tfrac{17}{\sqrt{10101}},\,\tfrac{7}{\sqrt{1677}}]$ \\
27 & $(\tfrac{7}{\sqrt{1677}},\ \tfrac{38}{\sqrt{49321}}]$ & $0.170935$ & $0.171107$ & 87 & $0.9552$ & $(\tfrac{7}{\sqrt{1677}},\,\tfrac{38}{\sqrt{49321}}]$ \\
28 & $(\tfrac{38}{\sqrt{49321}},\ \tfrac{17}{\sqrt{9709}}]$ & $0.171107$ & $0.172529$ & 86 & $0.9600$ & $(\tfrac{38}{\sqrt{49321}},\,\tfrac{17}{\sqrt{9709}}]$ \\
29 & $(\tfrac{17}{\sqrt{9709}},\ \tfrac{40}{\sqrt{53599}}]$ & $0.172529$ & $0.172775$ & 85 & $0.9515$ & $(\tfrac{17}{\sqrt{9709}},\,\tfrac{40}{\sqrt{53599}}]$ \\
30 & $(\tfrac{40}{\sqrt{53599}},\ \tfrac{1}{\sqrt{31}}]$ & $0.172775$ & $0.179605$ & 84 & $1.0161$ & $(\tfrac{40}{\sqrt{53599}},\,\tfrac{2}{3\sqrt{14}})$ \\
31 & $(\tfrac{1}{\sqrt{31}},\ \tfrac{14}{\sqrt{5917}}]$ & $0.179605$ & $0.182002$ & 78 & $0.9689$ & $(\tfrac{1}{\sqrt{31}},\,\tfrac{14}{\sqrt{5917}}]$ \\
32 & $(\tfrac{14}{\sqrt{5917}},\ \tfrac{89}{\sqrt{237553}}]$ & $0.182002$ & $0.182604$ & 77 & $0.9628$ & $(\tfrac{14}{\sqrt{5917}},\,\tfrac{89}{\sqrt{237553}}]$ \\
33 & $(\tfrac{89}{\sqrt{237553}},\ \tfrac{11}{\sqrt{3589}}]$ & $0.182604$ & $0.183614$ & 76 & $0.9609$ & $(\tfrac{89}{\sqrt{237553}},\,\tfrac{11}{\sqrt{3589}}]$ \\
34 & $(\tfrac{11}{\sqrt{3589}},\ \tfrac{7}{\sqrt{1443}}]$ & $0.183614$ & $0.184274$ & 75 & $0.9550$ & $(\tfrac{11}{\sqrt{3589}},\,\tfrac{7}{\sqrt{1443}}]$ \\
35 & $(\tfrac{7}{\sqrt{1443}},\ \tfrac{51}{\sqrt{76323}}]$ & $0.184274$ & $0.184605$ & 74 & $0.9457$ & $(\tfrac{7}{\sqrt{1443}},\,\tfrac{51}{\sqrt{76323}}]$ \\
36 & $(\tfrac{51}{\sqrt{76323}},\ \tfrac{1}{2\sqrt{7}}]$ & $0.184605$ & $0.188982$ & 73 & $0.9777$ & $(\tfrac{51}{\sqrt{76323}},\,\tfrac{1}{2\sqrt{7}}]$ \\
37 & $(\tfrac{1}{2\sqrt{7}},\ \tfrac{11}{\sqrt{3367}}]$ & $0.188982$ & $0.189571$ & 71 & $0.9568$ & $(\tfrac{1}{2\sqrt{7}},\,\tfrac{11}{\sqrt{3367}}]$ \\
38 & $(\tfrac{11}{\sqrt{3367}},\ \tfrac{2}{\sqrt{109}}]$ & $0.189571$ & $0.191565$ & 70 & $0.9633$ & $(\tfrac{11}{\sqrt{3367}},\,\tfrac{2}{\sqrt{109}}]$ \\
39 & $(\tfrac{2}{\sqrt{109}},\ \tfrac{7}{\sqrt{1273}}]$ & $0.191565$ & $0.196193$ & 68 & $0.9815$ & $(\tfrac{2}{\sqrt{109}},\,\tfrac{7}{\sqrt{1273}}]$ \\
40 & $(\tfrac{7}{\sqrt{1273}},\ \tfrac{59}{\sqrt{89713}}]$ & $0.196193$ & $0.196981$ & 67 & $0.9749$ & $(\tfrac{7}{\sqrt{1273}},\,\tfrac{59}{\sqrt{89713}}]$ \\
41 & $(\tfrac{59}{\sqrt{89713}},\ \tfrac{5}{7\sqrt{13}}]$ & $0.196981$ & $0.198107$ & 66 & $0.9714$ & $(\tfrac{59}{\sqrt{89713}},\,\tfrac{5}{7\sqrt{13}}]$ \\
42 & $(\tfrac{5}{7\sqrt{13}},\ \tfrac{71}{\sqrt{126697}}]$ & $0.198107$ & $0.199469$ & 65 & $0.9698$ & $(\tfrac{5}{7\sqrt{13}},\,\tfrac{71}{\sqrt{126697}}]$ \\
43 & $(\tfrac{71}{\sqrt{126697}},\ \tfrac{17}{\sqrt{7189}}]$ & $0.199469$ & $0.200500$ & 64 & $0.9648$ & $(\tfrac{71}{\sqrt{126697}},\,\tfrac{17}{\sqrt{7189}}]$ \\
44 & $(\tfrac{17}{\sqrt{7189}},\ \tfrac{19}{\sqrt{8827}}]$ & $0.200500$ & $0.202231$ & 63 & $0.9662$ & $(\tfrac{17}{\sqrt{7189}},\,\tfrac{19}{\sqrt{8827}}]$ \\
45 & $(\tfrac{19}{\sqrt{8827}},\ \tfrac{3}{\sqrt{217}}]$ & $0.202231$ & $0.203653$ & 62 & $0.9643$ & $(\tfrac{19}{\sqrt{8827}},\,\tfrac{3}{\sqrt{217}}]$ \\
46 & $(\tfrac{3}{\sqrt{217}},\ \tfrac{16}{\sqrt{6097}}]$ & $0.203653$ & $0.204909$ & 61 & $0.9605$ & $(\tfrac{3}{\sqrt{217}},\,\tfrac{16}{\sqrt{6097}}]$ \\
47 & $(\tfrac{16}{\sqrt{6097}},\ \tfrac{5}{\sqrt{589}}]$ & $0.204909$ & $0.206021$ & 60 & $0.9550$ & $(\tfrac{16}{\sqrt{6097}},\,\tfrac{5}{\sqrt{589}}]$ \\
48 & $(\tfrac{5}{\sqrt{589}},\ \tfrac{19}{91}]$ & $0.206021$ & $0.208791$ & 59 & $0.9645$ & $(\tfrac{5}{\sqrt{589}},\,\tfrac{19}{91}]$ \\
49 & $(\tfrac{19}{91},\ \tfrac{1}{\sqrt{21}}]$ & $0.208791$ & $0.218218$ & 58 & $1.0357$ & $(\tfrac{19}{91},\,\tfrac{2}{\sqrt{87}})$ \\
50 & $(\tfrac{1}{\sqrt{21}},\ \tfrac{1}{\sqrt{19}}]$ & $0.218218$ & $0.229416$ & 54 & $1.0658$ & $(\tfrac{1}{\sqrt{21}},\,\tfrac{2}{9})$ \\
51 & $(\tfrac{1}{\sqrt{19}},\ \tfrac{6}{\sqrt{651}}]$ & $0.229416$ & $0.235159$ & 48 & $0.9954$ & $(\tfrac{1}{\sqrt{19}},\,\tfrac{6}{\sqrt{651}}]$ \\
52 & $(\tfrac{6}{\sqrt{651}},\ \tfrac{17}{\sqrt{5187}}]$ & $0.235159$ & $0.236043$ & 47 & $0.9820$ & $(\tfrac{6}{\sqrt{651}},\,\tfrac{17}{\sqrt{5187}}]$ \\
53 & $(\tfrac{17}{\sqrt{5187}},\ \tfrac{3}{2\sqrt{39}}]$ & $0.236043$ & $0.240192$ & 46 & $0.9952$ & $(\tfrac{17}{\sqrt{5187}},\,\tfrac{3}{2\sqrt{39}}]$ \\
54 & $(\tfrac{3}{2\sqrt{39}},\ \tfrac{25}{\sqrt{10621}}]$ & $0.240192$ & $0.242581$ & 45 & $0.9930$ & $(\tfrac{3}{2\sqrt{39}},\,\tfrac{25}{\sqrt{10621}}]$ \\
55 & $(\tfrac{25}{\sqrt{10621}},\ \tfrac{7}{\sqrt{817}}]$ & $0.242581$ & $0.244899$ & 44 & $0.9896$ & $(\tfrac{25}{\sqrt{10621}},\,\tfrac{7}{\sqrt{817}}]$ \\
56 & $(\tfrac{7}{\sqrt{817}},\ \tfrac{7}{\sqrt{793}}]$ & $0.244899$ & $0.248577$ & 43 & $0.9964$ & $(\tfrac{7}{\sqrt{817}},\,\tfrac{7}{\sqrt{793}}]$ \\
57 & $(\tfrac{7}{\sqrt{793}},\ \tfrac{1}{4}]$ & $0.248577$ & $0.250000$ & 42 & $0.9844$ & $(\tfrac{7}{\sqrt{793}},\,\tfrac{1}{4}]$ \\
58 & $(\tfrac{1}{4},\ \tfrac{2}{3\sqrt{7}}]$ & $0.250000$ & $0.251976$ & 41 & $0.9762$ & $(\tfrac{1}{4},\,\tfrac{2}{3\sqrt{7}}]$ \\
59 & $(\tfrac{2}{3\sqrt{7}},\ \tfrac{3}{\sqrt{133}}]$ & $0.251976$ & $0.260133$ & 39 & $0.9897$ & $(\tfrac{2}{3\sqrt{7}},\,\tfrac{3}{\sqrt{133}}]$ \\
60 & $(\tfrac{3}{\sqrt{133}},\ \tfrac{5}{19}]$ & $0.260133$ & $0.263158$ & 38 & $0.9868$ & $(\tfrac{3}{\sqrt{133}},\,\tfrac{5}{19}]$ \\
61 & $(\tfrac{5}{19},\ \tfrac{1}{\sqrt{13}}]$ & $0.263158$ & $0.277350$ & 37 & $1.0673$ & $(\tfrac{5}{19},\,\tfrac{4}{\sqrt{222}})$ \\
62 & $(\tfrac{1}{\sqrt{13}},\ \tfrac{8}{\sqrt{817}}]$ & $0.277350$ & $0.279885$ & 35 & $1.0282$ & --- \\
63 & $(\tfrac{8}{\sqrt{817}},\ \tfrac{2}{7}]$ & $0.279885$ & $0.285714$ & 34 & $1.0408$ & $(\tfrac{8}{\sqrt{817}},\,\tfrac{2}{\sqrt{51}})$ \\
64 & $(\tfrac{2}{7},\ \tfrac{17}{\sqrt{3367}}]$ & $0.285714$ & $0.292973$ & 31 & $0.9978$ & $(\tfrac{2}{7},\,\tfrac{17}{\sqrt{3367}}]$ \\
65 & $(\tfrac{17}{\sqrt{3367}},\ \tfrac{29}{\sqrt{9331}}]$ & $0.292973$ & $0.300216$ & 30 & $1.0140$ & $(\tfrac{17}{\sqrt{3367}},\,\tfrac{2}{3\sqrt{5}})$ \\
66 & $(\tfrac{29}{\sqrt{9331}},\ \tfrac{5}{\sqrt{273}}]$ & $0.300216$ & $0.302614$ & 29 & $0.9959$ & $(\tfrac{29}{\sqrt{9331}},\,\tfrac{5}{\sqrt{273}}]$ \\
67 & $(\tfrac{5}{\sqrt{273}},\ \tfrac{3}{\sqrt{93}}]$ & $0.302614$ & $0.311086$ & 28 & $1.0161$ & $(\tfrac{5}{\sqrt{273}},\,\tfrac{2}{\sqrt{42}})$ \\
68 & $(\tfrac{3}{\sqrt{93}},\ \tfrac{12}{\sqrt{1443}}]$ & $0.311086$ & $0.315899$ & 27 & $1.0104$ & $(\tfrac{3}{\sqrt{93}},\,\tfrac{4}{9\sqrt{2}})$ \\
69 & $(\tfrac{12}{\sqrt{1443}},\ \tfrac{2}{\sqrt{39}}]$ & $0.315899$ & $0.320256$ & 26 & $1.0000$ & $(\tfrac{12}{\sqrt{1443}},\,\tfrac{2}{\sqrt{39}})$ \\
70 & $(\tfrac{2}{\sqrt{39}},\ \tfrac{1}{3}]$ & $0.320256$ & $0.333333$ & 25 & $1.0417$ & $(\tfrac{2}{\sqrt{39}},\,\tfrac{4}{5\sqrt{6}})$ \\
71 & $(\tfrac{1}{3},\ \tfrac{1}{\sqrt{7}}]$ & $0.333333$ & $0.377964$ & 24 & $1.2857$ & --- \\
72 & $(\tfrac{1}{\sqrt{7}},\ \tfrac{3}{\sqrt{57}}]$ & $0.377964$ & $0.397360$ & 18 & $1.0658$ & $(\tfrac{1}{\sqrt{7}},\,\tfrac{2}{3\sqrt{3}})$ \\
73 & $(\tfrac{3}{\sqrt{57}},\ \tfrac{17}{\sqrt{1729}}]$ & $0.397360$ & $0.408838$ & 17 & $1.0656$ & --- \\
74 & $(\tfrac{17}{\sqrt{1729}},\ \tfrac{8}{19}]$ & $0.408838$ & $0.421053$ & 16 & $1.0637$ & --- \\
75 & $(\tfrac{8}{19},\ \tfrac{2}{\sqrt{21}}]$ & $0.421053$ & $0.436436$ & 15 & $1.0714$ & $(\tfrac{8}{19},\,\tfrac{4}{3\sqrt{10}})$ \\
76 & $(\tfrac{2}{\sqrt{21}},\ \tfrac{2}{\sqrt{19}}]$ & $0.436436$ & $0.458831$ & 14 & $1.1053$ & --- \\
77 & $(\tfrac{2}{\sqrt{19}},\ \tfrac{1}{2}]$ & $0.458831$ & $0.500000$ & 13 & $1.2188$ & --- \\
78 & $(\tfrac{1}{2},\ \tfrac{5}{\sqrt{91}}]$ & $0.500000$ & $0.524142$ & 12 & $1.2363$ & --- \\
79 & $(\tfrac{5}{\sqrt{91}},\ \tfrac{7}{13}]$ & $0.524142$ & $0.538462$ & 11 & $1.1960$ & --- \\
80 & $(\tfrac{7}{13},\ \tfrac{1}{\sqrt{3}}]$ & $0.538462$ & $0.577350$ & 10 & $1.2500$ & --- \\
81 & $(\tfrac{1}{\sqrt{3}},\ \tfrac{3}{\sqrt{21}}]$ & $0.577350$ & $0.654654$ & 8 & $1.2857$ & --- \\
82 & $(\tfrac{3}{\sqrt{21}},\ \tfrac{2}{3}]$ & $0.654654$ & $0.666667$ & 7 & $1.1667$ & --- \\
83 & $(\tfrac{2}{3},\ 1]$ & $0.666667$ & $1.000000$ & 6 & $2.2500$ & --- \\
84 & $(1,\ \tfrac{2}{\sqrt{3}}]$ & $1.000000$ & $1.154701$ & 3 & $1.5000$ & --- \\
85 & $(\tfrac{2}{\sqrt{3}},\ 2]$ & $1.154701$ & $2.000000$ & 2 & $3.0000$ & --- \\
\end{longtable}
}

\end{document}